\documentclass[leqno,oneeqnum,onefignum,onetabnum]{siamltex1213}

\usepackage{amsmath,amssymb,mathtools}
\usepackage[sc]{mathpazo}
\usepackage{algorithm,algorithmic}
\usepackage[sort,compress]{cite}
\usepackage{graphicx}
\usepackage{multirow}
\usepackage{makecell}
\usepackage{textcomp}
\usepackage{rotating}
\usepackage{siunitx}
\usepackage{paralist}
\usepackage{booktabs}
\usepackage{lscape}
\usepackage[table]{xcolor}
\usepackage{geometry}
\usepackage[hyperpageref]{backref}

\newtheorem{remark}{Remark}

\graphicspath{{figures/}{tables/}}

\makeatletter
\@ifclassloaded{standalone}{}{
\hypersetup{
pdfauthor={Andreas Mang and George Biros},
linkcolor=red,
citecolor=orange}}
\makeatother%

\newcounter{runidnum}
\newcommand{\runid}{\color{orange}\stepcounter{runidnum}\#\therunidnum}
\newcommand{\resetrunid}{\setcounter{runidnum}{0}}

\newcommand{\figref}[1]{Fig.~\ref{#1}}
\newcommand{\tabref}[1]{Tab.~\ref{#1}}
\newcommand{\secref}[1]{\S\ref{#1}}

\newcommand{\runref}[1]{run~\##1}
\newcolumntype{R}{>{\columncolor{gray!20}}r}
\newcolumntype{L}{>{\columncolor{gray!20}}l}
\newcommand{\mcol}[2]{\multicolumn{#1}{r}{#2}}

\newcommand{\mrowrot}[2]{
\parbox[t]{2mm}{\multirow{#1}{*}{\rotatebox[origin=c]{90}{#2}}}
}

\newcommand{\F}[1]{\ensuremath{\mathcal{#1}}}    
\newcommand{\D}[1]{\ensuremath{\mathcal{#1}}}    
\newcommand{\ns}[1]{\ensuremath{\mathbf{#1}}}
\newcommand{\fs}[1]{\ensuremath{\mathcal{#1}}}
\newcommand{\idiv}{\ensuremath{\nabla\cdot}}
\newcommand{\igrad}{\ensuremath{\nabla}}
\newcommand{\ilap}{\rotatebox[origin=c]{180}{$\nabla$}}

\newcommand{\half}[1]{\frac{#1}{2}}
\newcommand{\ld}{\ensuremath{\mathrm{d}_t}} 
\renewcommand{\d}[1]{\mathop{}\!\mathrm{d}#1}
\newcommand{\dt}{\d{t}}

\newcommand{\bigO}{\mathcal{O}}
\newcommand{\p} {\partial}

\newcommand{\vect}[1]{\boldsymbol{#1}} 
\newcommand{\mat}[1]{\boldsymbol{#1}}  
\newcommand{\bipa}{\begin{inparaenum}[(\itshape i\upshape)]}
\newcommand{\eipa}{\end{inparaenum}}
\newcommand{\bipasub}{\begin{inparaenum}[(\itshape a\upshape)]}
\newcommand{\eipasub}{\end{inparaenum}}

\newcommand{\ipoint}[1]{\textit{\textbf{#1}}:}

\newcommand{\iut}[1]{\int_0^1#1\dt}

\newcommand{\defeq}{\ensuremath{\mathrel{\mathop:}=}}

\newcommand{\T}{\ensuremath{\mathsf{T}}}

\def\stitle{A NEWTON--KRYLOV SOLVER FOR CONSTRAINED DIFFEOMORPHIC IMAGE REGISTRATION}

\title{A SEMI-LAGRANGIAN TWO-LEVEL PRECONDITIONED NEWTON--KRYLOV SOLVER FOR CONSTRAINED DIFFEOMORPHIC IMAGE REGISTRATION\thanks{\textbf{Funding:} This material is based upon work supported by AFOSR grants FA9550-12-10484 and FA9550-11-10339; by NSF grant CCF-1337393; by the U.S. Department of Energy, Office of Science, Office of Advanced Scientific Computing Research, Applied Mathematics program under Award Numbers DE-SC0010518 and DE-SC0009286; by NIH grant 10042242; by DARPA grant W911NF-115-2-0121; and by the Technische Universit\"{a}t M\"{u}nchen---Institute for Advanced Study, funded by the German Excellence Initiative (and the European Union Seventh Framework Programme under grant agreement 291763). Any opinions, findings, and conclusions or recommendations expressed herein are those of the authors and do not necessarily reflect the views of the AFOSR, the DOE, the NIH, the DARPA, or the NSF.}}

\author{Andreas Mang\thanks{The Institute for Computational Engineering and Sciences, The University of Texas at Austin, Austin, Texas, 78712, US. Current address: Department of Mathematics at the University of Houston, Houston, Texas 77004, US; {\tt andreas@math.uh.edu}.} \and George Biros\thanks{The Institute for Computational Engineering and Sciences, The University of Texas at Austin, Austin, Texas, 78712, US; {\tt gbiros@acm.org}}}

\begin{document}

\maketitle

\begin{abstract}
We propose an efficient numerical algorithm for the solution of diffeomorphic image registration problems. We use a variational formulation constrained by a partial differential equation (PDE), where the constraints are a scalar transport equation.

We use a pseudospectral discretization in space and second-order accurate semi-Lagrangian time stepping scheme for the transport equations. We solve for a stationary velocity field using a preconditioned, globalized, matrix-free Newton-Krylov scheme. We propose and test a two-level Hessian preconditioner. We consider two strategies for inverting the preconditioner on the coarse grid: a nested preconditioned conjugate gradient method (exact solve) and a nested Chebyshev iterative method (inexact solve) with a fixed number of iterations.

We test the performance of our solver in different synthetic and real-world two-dimensional application scenarios. We study grid convergence and computational efficiency of our new scheme. We compare the performance of our solver against our initial implementation that uses the same spatial discretization but a standard, explicit, second-order Runge-Kutta scheme for the numerical time integration of the transport equations and a single-level preconditioner. Our improved scheme delivers significant speedups over our original implementation. As a highlight, we observe a 20$\times$ speedup for a two dimensional, real world multi-subject medical image registration problem.
\end{abstract}

\newcommand{\slugmaster}{\slugger{siims}{xxxx}{xx}{x}{x--x}} 



\pagestyle{myheadings}
\thispagestyle{plain}
\markboth
{
ANDREAS MANG AND GEORGE BIROS
}
{
\stitle
}

\section{Introduction}

Image registration finds numerous applications in image analysis and computer vision~\cite{Hajnal:2001a,Sotiras:2013a}. Image registration establishes \emph{meaningful} spatial correspondence between two images $m_R:\bar{\Omega}\rightarrow\ns{R}$ (the `'reference image`') and $m_T:\bar{\Omega}\rightarrow\ns{R}$ (the `'template image`') of a scene such that the deformed template image $m_T$ becomes similar to $m_R$, i.e., $m_T\circ\vect{y}\approx m_R$~\cite{Modersitzki:2004a}; the images are defined on an open set $\Omega\subset\ns{R}^d$, $d\in\{2,3\}$, with closure $\bar{\Omega}\defeq \Omega\cup\partial\Omega$ and boundary $\partial\Omega$, $\circ$ denotes the composition of two functions, and $\vect{y}:\bar{\Omega}\rightarrow\bar{\Omega}$ is the sought after deformation map. There exist various approaches to image registration; we refer to~\cite{Modersitzki:2004a,Fischer:2008a,Sotiras:2013a} for a lucid overview.

Image registration is typically formulated as a variational optimization problem with an objective functional that consists of a data fidelity term and a Tikhonov-type regularization norm~\cite{Amit:1994a}; the unregularized problem is ill-posed. Here, we follow up on our preceding work on \emph{constrained diffeomorphic image registration}~\cite{Mang:2016a,Mang:2015a}. In diffeomorphic image registration we require that the map $\vect{y}$ is a \emph{diffeomorphism}, i.e., $\vect{y}$ is a bijection, continuously differentiable, and has a continuously differentiable inverse. Formally, we require that $\det\igrad\vect{y} \not=0, \igrad\vect{y}\in\ns{R}^{d\times d}, \,\forall \vect{x}\in\Omega$, and---under the assumption that $\vect{y}$ is orientation preserving---$\det\igrad\vect{y} > 0$, $\forall \vect{x}\in\Omega$.

Different approaches to guarantee a diffeomorphic $\vect{y}$ have appeared in the past. One approach is to penalize $\det\igrad\vect{y}$ as done in~\cite{Burger:2013a,Droske:2003a,Haber:2004a,Haber:2007a,Sdika:2008a}. Another approach is to change the formulation; instead of inverting directly for the deformation map $\vect{y}$, we invert for its velocity $\vect{v}=\d_t \vect{y}$. If $\vect{v}$ is sufficiently smooth it can be guaranteed that the resulting $\vect{y}$ is a diffeomorphism~\cite{Beg:2005a,Dupuis:1998a,Trouve:1998a}. In our formulation, we augment this type of smoothness regularization by constraints on the divergence of $\vect{v}$~\cite{Mang:2016a,Mang:2015a}. For instance, for $\idiv\vect{v}=0$ the flow becomes incompressible. This is equivalent to enforcing $\det\igrad\vect{y} = 1$~\cite[pages~77ff.]{Gurtin:1981a}.

Velocity field formulations for diffeomorphic image registration can be distinguished between approaches that invert for a \emph{time dependent} $\vect{v}$~\cite{Beg:2005a,Dupuis:1998a,Hart:2009a,Mang:2015a} and approaches that invert for a \emph{stationary} $\vect{v}$~\cite{Hernandez:2009a,Lorenzi:2013a,Mang:2016a}. We invert for a stationary $\vect{v}$. We formulate the diffeomorphic image registration problem as a PDE constrained optimization problem, where the constraint is a transport equation for the scalar field $m : \bar{\Omega}\times[0,1]\rightarrow\ns{R}$ (the image intensities). Due to ill-conditioning, non-convexity, large-problem size, infinite-dimensional structure, and the need for adjoint operators, such problems are challenging to solve.  We use a reduced space Newton--Krylov method~\cite{Mang:2015a}. In reduced space methods we eliminate state variables (in our case the transported image) and iterate in the control variable space (in our case the velocity space). Newton methods typically display faster convergence than gradient-descent methods (see~\cite{Mang:2015a}). Using a Newton method, however, requires solving linear systems with the reduced space Hessian, which---upon discretization---is a large, dense, and ill-conditioned operator. Efficient preconditioning is critical for making our solver effective across a wide spectrum of image resolutions, regularization weights, and inversion tolerances. Standard preconditioning techniques like incomplete factorization cannot be applied since we do not have access to the matrix entries (too expensive to compute). Instead, we present a matrix-free, two-level preconditioner for the reduced space Hessian that significantly improves performance. Another computational challenge of our formulation is that the reduced space formulation requires the exact solution of two hyperbolic transport equations---the state and adjoint equations of our problem---every time we evaluate the reduced gradient or apply the reduced space Hessian operator. We introduce a semi-Lagrangian formulation to further speed up our solver.

\subsection{Outline of the Method}
\label{s:outline}

We are given two functions $m_R : \bar{\Omega}\rightarrow\ns{R}$ (fixed image) and $m_T : \bar{\Omega}\rightarrow\ns{R}$ (deformable image) compactly supported on an open set $\Omega \defeq (-\pi,\pi)^d$, $d\in\{2,3\}$, with boundary $\p\Omega$, and closure $\bar{\Omega}\defeq \Omega \cup \p\Omega$. We solve for a stationary velocity field $\vect{v}\in\fs{U}$ and a mass source $w\in\fs{W}$ as follows~\cite{Mang:2016a}:
\begin{subequations}
\label{e:ip}
\begin{equation}
\min_{m,\vect{v},w}
\half{1}
\|m_R - m_1\|^2_{L^2(\Omega)}
+ \half{\beta_v} \|\vect{v}\|^2_{\fs{V}}
+ \half{\beta_w} \|w\|^2_{\fs{W}}
\label{e:ip:objective}
\end{equation}
\noindent subject to
\begin{align}
\p_t m + \igrad m \cdot \vect{v} & = 0
&& {\rm in}\;\Omega\times(0,1],
\label{e:ip:transport}
\\
m & = m_T
&& {\rm in}\;\Omega\times\{0\},
\label{e:ip:transport-ic}
\\
\idiv \vect{v} &= w
&& {\rm in} \; \Omega,
\label{e:ip:constraint}
\end{align}
\end{subequations}

\noindent and periodic boundary conditions on $\partial\Omega$. In our formulation $m_1(\vect{x})\defeq m(\vect{x},t=1)$----i.e., the solution of the hyperbolic transport equation~\eqref{e:ip:transport} with initial condition~\eqref{e:ip:transport-ic}---is equivalent to $m_T\circ\vect{y}$; the deformation map $\vect{y}$ can be computed from $\vect{v}$ in a post-processing step (see, e.g.,~\cite{Mang:2015a,Mang:2016a}). The weights $\beta_v>0$, and $\beta_w>0$ control the regularity of $\vect{v}$.

The regularization norm for $\vect{v}$ not only alleviates issues related to the ill-posedness of our problem but also ensures the existence of a diffeomorphism $\vect{y}$ parameterized by $\vect{v}$ if chosen appropriately. The constraint in~\eqref{e:ip:constraint} allows us to control volume change; setting $w=0$ results in an incompressible diffeomorphism $\vect{y}$, i.e., the deformation gradient $\det\igrad\vect{y}$ is fixed to one for all $\vect{x}\in\Omega$. The deformation map $\vect{y}$ is no longer incompressible if we allow $w$ to deviate from zero (this formulation has originally been introduced in~\cite{Mang:2016a}; a similar formulation can be found in \cite{Borzi:2002a}). We can control this deviation with $\beta_w$; the regularization norm for $w$ acts like a penalty on $\idiv \vect{v}$. We will specify and discuss the choices for the spaces $\fs{U}$, $\fs{V}$, and $\fs{W}$ in more detail~\secref{s:formulation} and~\secref{s:optsys-details}.

We use the method of Lagrange multipliers to solve~\eqref{e:ip}. We first formally derive the optimality conditions and then discretize using a pseudospectral discretization in space with a Fourier basis (i.e., we use an \emph{optimize-then-discretize} approach; see \secref{s:solver}). We solve for the first-order optimality conditions using a globalized, matrix-free, preconditioned, inexact Newton--Krylov algorithm for the velocity field $\vect{v}$ (see~\cite{Mang:2015a} for details). The hyperbolic transport equations are solved via a semi-Lagrangian method.

\subsection{Contributions}

Our Newton-Krylov scheme has originally been described in~\cite{Mang:2015a}, in which we compared it to gradient-descent approach in the Sobolev space induced by the regularization operator (the latter approach is, e.g., used in~\cite{Hart:2009a}). The latter, as expected, is extremely slow and not competitive with (Gauss--)Newton schemes. In~\cite{Mang:2016a} we introduced and studied different regularization functionals, and compared the performance of our method against existing approaches for diffeomorphic image registration, in particular the Demons family of algorithms~\cite{Vercauteren:2008a,Vercauteren:2009a}. Here we extend our preceding work in the following ways:
\begin{itemize}
\item We propose a semi-Lagrangian formulation for our entire optimality system, i.e., the state, adjoint, and incremental state and adjoint equations. We compare it with an stabilized Runge--Kutta method, which we also introduce here; we show that the semi-Lagrangian scheme has excellent stability properties.
\item We introduce an improved preconditioner for the reduced Hessian system. It is a two-level preconditioner that uses spectral restriction and prolongation operators and a Chebyshev stationary iterative method for an approximate coarse grid solve.
\item We provide an experimental study of the performance of our improved numerical scheme based on synthetic and real-world problems. We study self-convergence, grid convergence, numerical accuracy, convergence as a function of the regularization parameter, and the time to solution. We account for different constraints and regularization norms.
\end{itemize}

Taken together, the new algorithm results in order of magnitude speedups over the state-of-the-art. For example, for a magnetic resonance image of a brain with $512^2$ resolution the new scheme is 18$\times$ faster (see~\tabref{t:preconditioners-kktsolve-convergence} in \secref{s:experiments}) than the scheme described in~\cite{Mang:2016a}.

\subsection{Limitations and Unresolved Issues}

Several limitations and unresolved issues remain. We assume similar intensity statistics for $m_R$ and $m_T$. This is a common assumption in many deformable registration algorithms. For multimodal registration problems we have to replace the squared $L^2$-distance in~\eqref{e:ip:objective} with more involved distance measure; examples can be found in~\cite{Modersitzki:2004a,Sotiras:2013a}. We present results only for $d=2$. Nothing in our formulation and numerical approximation is specific to the two-dimensional case. In this work we discuss improvements of the algorithm used in our preceding work~\cite{Mang:2016a,Mang:2015a} en route to an effective three-dimensional solver. Once this three-dimensional solver is available, we will extend the study presented in~\cite{Mang:2015a}, by providing a detailed comparison of our method against diffeomorphic image registration approaches of other groups in terms of efficiency and inversion accuracy.

\subsection{Related Work}

The body of literature on diffeomorphic image registration, numerical optimization in optimal control, preconditioning of KKT systems, and the effective solution of hyperbolic transport equations is extensive. We limit the discussion to work that is most relevant to ours.

\subsubsection{Diffeomorphic Image Registration}

Lucid overviews for image registration can be found in~\cite{Fischer:2008a,Modersitzki:2009a,Sotiras:2013a}. Related work on velocity field based diffeomorphic image registration is discussed in~\cite{Ashburner:2007a,Ashburner:2011a,Beg:2005a,Borzi:2002a,Hart:2009a,Lee:2010a,Mang:2016a,Mang:2015a} and references therein. Related optimal control formulations for image registration are described in~\cite{Barbu:2016a,Benzi:2011a,Borzi:2002a,Chen:2011a,Lee:2010a,Mang:2016a,Mang:2015a,Simoncini:2012a,Vialard:2012a,Ruhnau:2007a}. Most work on velocity based diffeomorphic registration considers first order information for numerical optimization (see, e.g., \cite{Beg:2005a,Borzi:2002a,Cao:2005a,Chen:2011a,Hart:2009a,Lee:2010a,Vialard:2012a}), with the exceptions of our own work~\cite{Mang:2016a,Mang:2015a} and~\cite{Ashburner:2011a,Benzi:2011a,Hernandez:2014a,Simoncini:2012a}; only~\cite{Benzi:2011a,Simoncini:2012a} discuss preconditioning strategies (see also below). The application of a Newton--Krylov solver for incompressible and near-incompressible formulations (with an additional control on a mass-source term) for diffeomorphic image registration is, to the best of our knowledge, exclusive to our group~\cite{Mang:2016a,Mang:2015a}.

\subsubsection{PDE Constrained Optimization}

There exists a huge body of literature for the numerical solution of PDE constrained optimization problems. The numerical implementation of an efficient solver is, in many cases, tailored towards the nature of the control problem, e.g., by accounting for the type and structure of the PDE constraints; see for instance~\cite{Adavani:2008b,Biros:2008a} (elliptic),~\cite{Adavani:2008a,Mang:2012b,Gholami:2016a,Stoll:2015a} (parabolic), or~\cite{Benzi:2011a,Borzi:2002a,Lee:2010a} (hyperbolic). We refer to~\cite{Biegler:2003a,Borzi:2012a,Gunzburger:2003a,Herzog:2010a,Hinze:2009a} for an overview on theoretical and algorithmic developments in optimal control and PDE constrained optimization. A survey on strategies for preconditioning saddle point problems can be found in~\cite{Benzi:2005a}. We refer to~\cite{Borzi:2009a} for an overview on multigrid methods for optimal control problems.

Our preconditioner can be viewed as a simplified two level multigrid v-cycle with a smoother based on the inverse regularization operator and the coarse grid solve is inexact. We note, that more sophisticated multigrid preconditioners for the reduced Hessian exist~\cite{Adavani:2008a,Borzi:2002a}. Multigrid approaches have been considered in~\cite{Borzi:2002a} for optical flow and in~\cite{Benzi:2011a,Simoncini:2012a} for the Monge-Kantorovich functional. The work of~\cite{Borzi:2002a} is the most pertinent to our problem. It is a space-time multigrid in the full KKT conditions and the time discretization scheme is CFL restricted, and, thus, very expensive. The effectiveness of the smoother depends on the regularization functional—--it is unclear how to generalize it to incompressible velocities. Our scheme is simpler to implement, supports general regularizations, and is compatible with our semi-Lagrangian time discretization. The preconditioner in~\cite{Benzi:2011a} is a block triangular preconditioner based on a perturbed representation of the GN approximation of the full space KKT system. A similar preconditioner that operates on the reduced space Hessian is described in~\cite{Simoncini:2012a}. In some sense, we do not approximate the structure of our Hessian operator; we invert an exact representation. We amortize the associated costs as follows: \bipa\item we solve for the action of the inverse inexactly, and \item we invert the operator on a coarser grid\eipa.

\subsubsection{The Semi-Lagrangian Method}

We refer to~\cite{Ewing:2001a} for a summary on solvers for advection dominated systems. Example implementations for the solution of hyperbolic transport equations that have been considered in the work cited above are implicit Lax-Friedrich schemes~\cite{Benzi:2011a,Simoncini:2012a}, explicit high-order total variation diminishing schemes~\cite{Borzi:2002a,Chen:2011a,Hart:2009a}, or explicit, pseudospectral (in space) RK2 schemes~\cite{Mang:2016a,Mang:2015a}. These schemes suffer either from numerical diffusion and/or CFL time step restrictions. We use a high-order, unconditionally stable semi-Lagrangian formulation. Semi-Lagrangian methods are well established and have first been considered in numerical weather prediction~\cite{Staniforth:1991a}. The use of semi-Lagrangian schemes is not new in the context of diffeomorphic image registration. However, such schemes have only been used to solve for the deformation map and/or solve the forward problem~\cite{Beg:2005a,Cao:2005a,Chen:2011a,Hernandez:2009a} and for the adjoint problem in the context of approximate gradient-descent methods.

\subsection{Organization and Notation}

We summarize our notation in \tabref{t:notation}. We summarize the  optimal control formulation for diffeomorphic image registration in~\secref{s:formulation}. We describe the solver in~\secref{s:solver}. We provide the optimality system and the Newton step in~\secref{s:optsys}. We describe the discretization in~\secref{s:discretization}. The schemes for integrating the hyperbolic PDEs that appear in our formulation are discussed in~\secref{s:timeintegration}. We describe our Newton--Krylov solver in~\secref{s:newtonkrylovmethod}; this includes a discussion of the preconditioners for the solution of the reduced space KKT system. We provide numerical experiments in \secref{s:experiments}. We conclude with~\secref{s:conclusions}.

\begin{table}
\caption{Commonly used notation and symbols.}
\centering
\begin{footnotesize}
\begin{tabular}{ll}
\toprule
Symbol/Notation & Description
\\\midrule
CFL             & Courant--Friedrichs--Lewy (condition)\\
FFT             & Fast Fourier Transform\\
GN              & Gauss--Newton\\
KKT             & Karush--Kuhn--Tucker (system)\\
matvec          & (Hessian) matrix-vector product\\
PCG             & Preconditioned Conjugate Gradient (method)\\
PCG($\epsilon$) & PCG, where $\epsilon>0$ indicates the used tolerance\\
PDE             & partial differential equation\\
PDE solve       & solution of a hyperbolic transport equation\\
RK2             & 2nd order Runge--Kutta (method)\\
RK2($c$)        & RK2 method, where $c$ indicates the employed CFL number\\
RK2A            & RK2 scheme based on an antisymmetric form\\
RK2A($c$)       & RK2A method, where $c$ indicates the employed CFL number\\
SL              & semi-Lagrangian (method)\\
SL($c$)         & SL method, where $c$ indicates the employed CFL number\\
\midrule
$d$             & spatial dimensionality; typically $d\in\{2,3\}$\\
$\Omega$        & spatial domain; $\Omega\defeq(-\pi,\pi)^d\subset\ns{R}^d$ with boundary $\p\Omega$ and closure $\bar{\Omega}\defeq\Omega\cup\p\Omega$\\
$\vect{x}$      & spatial coordinate; $\vect{x}\defeq(x^1,\ldots,x^d)^\T\in\ns{R}^d$\\
$m_R$           & reference image; $m_R : \bar{\Omega} \rightarrow \ns{R}$\\
$m_T$           & template image; $m_T : \bar{\Omega} \rightarrow \ns{R}$\\
$m$             & state variable (transported intensities); $m : \bar{\Omega}\times[0,1]\rightarrow\ns{R}$\\
$m_1$           & deformed template image (state variable at $t=1$); $m_1 : \bar{\Omega} \rightarrow \ns{R}$\\
$\lambda$       & adjoint variable (transport equation); $\lambda : \bar{\Omega}\times[0,1]\rightarrow\ns{R}$\\
$p$             & adjoint variable (incompressibility constraint); $p : \bar{\Omega}\rightarrow\ns{R}$\\
$\vect{v}$      & control variable (stationary velocity field); $\vect{v} : \bar{\Omega} \rightarrow \ns{R}^d$\\
$w$             & control variable (mass source); $w : \bar{\Omega} \rightarrow \ns{R}$\\
$\vect{b}$      & body force; $\vect{b} : \bar{\Omega} \rightarrow \ns{R}^d$\\
$\D{H}$         & (reduced) Hessian\\
$\vect{g}$      & (reduced) gradient\\
$\mat{y}$       & Eulerian (pullback) deformation map\\
$\mat{F}$       & deformation gradient at $t=1$ (computed from $\vect{v}$); $\mat{F} : \bar{\Omega}\rightarrow\ns{R}^{d\times d}$; $\mat{F}\defeq(\igrad\vect{y})^{-1}$\\
$\beta_v$       & regularization parameter for the control $\vect{v}$\\
$\beta_w$       & regularization parameter for the control $w$\\
$\D{A}$         & regularization operator (variation of regularization model acting on $\vect{v}$) \\
$\p_i$          & partial derivative with respect to $x^i$, $i=1,\ldots,d$\\
$\p_t$          & partial derivative with respect to time\\
$\ld$           & Lagrangian derivative\\
$\igrad$        & gradient operator (acts on scalar and vector fields)\\
$\ilap$         & Laplacian operator (acts on scalar and vector fields)\\
$\idiv$         & divergence operator (acts on vector and 2nd order tensor fields)\\
$\langle\cdot,\cdot\rangle_{L^2(\fs{X})}$ & $L^2$ inner product on $\fs{X}$
\\\bottomrule
\end{tabular}
\end{footnotesize}
\label{t:notation}
\end{table}

\section{Optimal Control Formulation}
\label{s:formulation}

We consider a PDE constrained formulation, where the constraints consist of a scalar transport equation for the image intensities. We solve for a stationary velocity field $\vect{v}\in\fs{U}$ and a mass-source $w\in\fs{W}$ as follows~\cite{Mang:2016a}:
\begin{subequations}
\label{e:varopt}
\begin{equation}
\label{e:varopt:objective}
\min_{m,\vect{v},w}
\fs{J}[\vect{v},w]=
\half{1}\|m_1 - m_R\|^2_{L^2(\Omega)}
+ \half{\beta_v}\|\vect{v}\|^2_{\fs{V}}
+ \half{\beta_w}\|w\|^2_{\fs{W}}
\end{equation}
\noindent subject to
\begin{align}
\p_t m + \igrad m\cdot\vect{v}
& = 0 && {\rm in}\;\Omega\times(0,1],
\label{e:varopt:transport}
\\
m & = m_T
&& {\rm in}\;\Omega\times\{0\},
\label{e:varopt:transport-ic}
\\
\idiv \vect{v} & = w
&& {\rm in} \; \Omega
\label{e:varopt:constraint}
\end{align}
\end{subequations}

\noindent and periodic boundary conditions on $\p\Omega$. We measure the similarity between the reference image $m_R$ and the deformed template image $m_1$ using a squared $L^2$-distance. The contributions of the regularization models for $w$ and $\vect{v}$ are controlled by the weights $\beta_v>0$ and $\beta_w>0$, respectively. We consider an $H^1$-regularization norm for $w$, i.e.,
\begin{equation}
\label{e:varopt:regw}
\|w\|^2_{H^1(\Omega)}
\defeq \int_{\Omega} \igrad w \cdot \igrad w + w^2 \d{\vect{x}}.
\end{equation}

\noindent We consider three quadratic regularization models for $\vect{v}$; an $H^1$-, an $H^2$-, and an $H^3$-seminorm:
\begin{equation}
\label{e:varopt:regv}
|\vect{v}|^2_{H^1(\Omega)^d}
\defeq \int_{\Omega} \igrad\vect{v} : \igrad\vect{v}\d{\vect{x}},
\quad
|\vect{v}|^2_{H^2(\Omega)^d}
\defeq\int_{\Omega} \ilap\vect{v} \cdot \ilap\vect{v}\d{\vect{x}},
\quad
\text{and}\quad
|\vect{v}|^2_{H^3(\Omega)^d}
\defeq \int_{\Omega} \igrad\ilap\vect{v} : \igrad\ilap\vect{v}\d{\vect{x}}.
\end{equation}

The use of an $H^1$-seminorm is motivated by related work in computational fluid dynamics; we will see that the first order variations of our formulation will result in a system that reflects a linear Stokes model under the assumption that we enforce $\idiv \vect{v} = 0$~\cite{Chen:2011a,Mang:2016a,Mang:2015a,Ruhnau:2007a}. We use an $H^2$-seminorm if we neglect the incompressibility constraint~\eqref{e:varopt:constraint}. This establishes a connection to related formulations for diffeomorphic image registration~\cite{Beg:2005a,Hart:2009a,Hernandez:2009a}; an $H^2$-norm is the paramount model in many algorithms (or its approximation via its Green's function; a Gaussian kernel)~\cite{Beg:2005a}.

\begin{remark}
The norm on $w$ acts like a penalty on $\idiv\vect{v}$. In fact, we can eliminate~\eqref{e:varopt:constraint} from~\eqref{e:varopt} by inserting $\idiv \vect{v}$ for $w$ into the regularization norm in~\eqref{e:varopt:objective}. If we neglect the incompressibility constraint~\eqref{e:varopt:constraint} the space $\fs{U}$ for $\vect{v}$ is given by the Sobolev space $\fs{V}$ (this formulation is, e.g., used in~\cite{Hart:2009a} for a non-stationary velocity with $H^2$-regularity in space and $L^2$ regularity in time). If we set $w$ in~\eqref{e:varopt:constraint} to zero, the computed velocity will be in the space of divergence free velocity fields with Sobolev regularity in space, as defined by $\fs{V}$ (examples for this formulation can be found in~\cite{Chen:2011a,Mang:2016a,Mang:2015a,Ruhnau:2007a}). For a non-zero $w$ we additionally require that the divergence of $\vect{v}$ is in $\fs{W}$. An equivalent formulation is, e.g., presented in~\cite{Borzi:2002a,Borzi:2002b}. They use $H^1$-regularity for $\vect{v}$ and stipulate $L^2$-regularity for its divergence, and proof existence of the state and adjoint variables for smooth images~\cite{Borzi:2002a}. In particular, they provide existence results for a unique, $H^1$-regular solution of the forward problem under the assumption of $H^1$-regularity for the template image. The same regularity requirements hold true for the adjoint equation. In our formulation, we not only require $\vect{v}$ to be an $H^1$-function, but also that its divergence is in $H^1$ (according to~\eqref{e:varopt:regw}). Another approach to impose regularity on $\vect{v}$ is to not only control the divergence but also control its curl (see, e.g.,~\cite{Lee:2010a,Amrouche:1998a}). We provide additional remarks in~\secref{s:optsys-details}.
\end{remark}

\section{Numerics and Solver}
\label{s:solver}

In what follows, we describe our numerical solver for computing a discrete approximation to the continuous problem. We use a globalized, preconditioned, inexact, reduced space\footnote{By \emph{reduced space} we mean that we will only iterate on the reduced space of the velocity $\vect{v}$; we assume that the state and adjoint equations are fulfilled exactly. This is different to \emph{all-at-once} or \emph{full space approaches}, in which one iterates on all unknown variables simultaneously (see~\secref{s:optsys} and~\secref{s:optsys-details}).} (Gauss--)Newton--Krylov method. Our scheme is described in detail in~\cite{Mang:2015a}. We will briefly recapitulate the key ideas and main building blocks.

We use the (formal) Lagrangian method~\cite{Lions:1971a} to solve~\eqref{e:varopt}; the Lagrangian functional $\D{L}$ is given by
\begin{align}
\label{e:lagrangian}
\F{L}[\vect{\phi}] \defeq
& \F{J}[\vect{v},w]
+ \int_0^1\langle\p_t m + \vect{v}\cdot\igrad m,\lambda\rangle_{L^2(\Omega)}\d{t}
+ \langle m_0 - m_T,\upsilon\rangle_{L^2(\Omega)}
- \langle \idiv \vect{v} - w,p\rangle_{L^2(\Omega)}
\end{align}

\noindent with $\vect{\phi}\defeq(m,\lambda,p,w,\vect{v})$ and Lagrange multipliers $\lambda:\bar{\Omega}\times[0,1]\rightarrow\ns{R}$ for the hyperbolic transport equation~\eqref{e:varopt:transport}, $\nu:\bar{\Omega}\rightarrow\ns{R}$ for the initial condition~\eqref{e:varopt:transport-ic}, and $p:\bar{\Omega}\rightarrow\ns{R}$ for the incompressibility constraint~\eqref{e:varopt:constraint} (we neglect the periodic boundary conditions for simplicity). The Lagrange multiplier functions inherit the boundary conditions of the forward operator.

\begin{remark}
We can consider two numerical strategies to tackle~\eqref{e:varopt}. We can either use an \emph{optimize-then-discretize} approach or a \emph{discretize-then-optimize} approach. We choose the former, i.e., we compute variations of the continuous problem and then discretize the optimality system. In general, this approach does not guarantee that the discretization of the gradient is consistent with the discretized objective. Further, it is not guaranteed that the discretized forward and adjoint operators are transposes of one another. Likewise, it is not guaranteed that the discretized Hessian is a symmetric operator. We report numerical experiments to quantify these errors; we will see that they are below the tolerances we target for the inversion. By using a discretize-then-optimize approach one can (by construction) guarantee that the derived operators are consistent. However, it is, e.g., not guaranteed that the forward and adjoint operators (in the transposed sense) yield the same numerical accuracy (see, e.g.,~\cite{Hager:2000a,Hager:2000b}). We refer, e.g., to~\cite{Gunzburger:2003a,Borzi:2012a} for additional remarks on the discretization of optimization and control problems.
\end{remark}

\subsection{Optimality Conditions and Newton Step}
\label{s:optsys}

From Lagrange multiplier theory we know that we require vanishing variations of $\fs{L}$ in~\eqref{e:lagrangian} with respect to the state, adjoint, and control variables $\vect{\phi}$ for an admissible solution to~\eqref{e:varopt}. We present the steps necessary to evaluate the reduced gradient and Hessian matvec. The associated PDE operators are derived using calculus of variations, and invoking Green's identities. We will see that the optimality conditions of our problem form a system of PDEs. This system needs to be solved to find a solution of~\eqref{e:varopt}. We will only present the strong form of our reduced space formulation.\footnote{We refer to~\cite{Biros:2005a,Biros:2005b} for more details on reduced-space methods.} Note, that we also eliminate the incompressibility constraint from the optimality system (see~\cite{Mang:2015a,Mang:2016a} for details; we comment on this in more detail in~\secref{s:optsys-details}); we only iterate on the reduced space for the velocity field $\vect{v}$. The expression for the reduced gradient for our problem is given by
\begin{equation}
\vect{g}(\vect{v}) \defeq
  \beta_v\D{A}[\vect{v}] + \D{K}[\vect{b}]
= \beta_v\D{A}[\vect{v}] + \D{K}[\iut{\lambda\igrad m}]
\label{e:control-elim}
\end{equation}

\noindent with (pseudo-)differential operators $\D{A}$ (regularization) and $\D{K}$ (projection); the definitions are given below. Formally, we require $\vect{g}(\vect{v}^\star) = \vect{0}$ for an admissible solution $\vect{v}^\star$ to~\eqref{e:varopt}. We can compute this minimizer iteratively using $\vect{g}$ in a gradient descent scheme. To evaluate $\vect{g}$ we need to find the space-time fields $m$ and $\lambda$ given a candidate $\vect{v}$. We can compute $m$ by solving the \emph{state equation} (primal)
\begin{subequations}
\label{e:state-elim}
\begin{align}
\p_t m +\vect{v}\cdot\igrad m & = 0
&&{\rm in} \;\; \Omega \times (0,1],
\label{e:state-pde-elim}
\\
m  &= m_T
&&{\rm in} \;\; \Omega \times\{0\},
\label{e:state-pde-ic-elim}
\end{align}
\end{subequations}

\noindent with periodic boundary conditions on $\p\Omega$ \emph{forward} in time. Once we have found $m$ at $t=1$ we can compute $\lambda$ by solving the \emph{adjoint} or \emph{costate equation} (dual)
\begin{subequations}
\label{e:adj-elim}
\begin{align}
-\p_t\lambda - \idiv \lambda\vect{v} & = 0
&&\text{in}\;\;\Omega\times[0,1),
\label{e:adj-pde-elim}
\\
\lambda &= -(m - m_R)
&&{\rm in} \;\; \Omega \times\{1\},
\label{e:adj-pde-fc-elim}
\end{align}
\end{subequations}

\noindent with periodic boundary conditions on $\p\Omega$ \emph{backward} in time; for vanishing $\idiv\vect{v}$ \eqref{e:adj-pde-elim} will also be a transport equation.

What is missing to complete the picture for $\vect{g}$ is a specification of the operators $\D{A}$ and $\D{K}$. The differential operator $\D{A}$ in~\eqref{e:control-elim} corresponds to the first variation of the seminorms in~\eqref{e:varopt:regv}. We have
\begin{equation}
\label{e:diffopregv}
\D{A}[\vect{v}] = -\ilap \vect{v},
\qquad
\D{A}[\vect{v}] = \ilap^2\vect{v},
\qquad
\text{and}
\qquad
\D{A}[\vect{v}] = \ilap^3\vect{v}
\end{equation}

\noindent for the $H^1$, $H^2$, and $H^3$ case, respectively, resulting in an \emph{elliptic}, \emph{biharmonic}, or \emph{triharmonic} integro-differential control equation for $\vect{v}$, respectively. The pseudo-differential operator $\D{K}$ in \eqref{e:control-elim} originates from the elimination of $p$ and~\eqref{e:ip:constraint}. For instance, if we set $w=0$ we obtain the Leray operator $\D{K}[\vect{b}] \defeq -\igrad\ilap^{-1}\idiv\vect{b} + \vect{b}$; for non-zero $w$ this operator becomes more complicated (see~\cite{Mang:2015a,Mang:2016a} for details on the derivation of this operator). Combining the state (primal)~\eqref{e:state-elim}, the adjoint~\eqref{e:adj-elim}, and the control equation~\eqref{e:control-elim} provides the formal optimality conditions (see \secref{s:optsys-details}).

A common strategy to compute a minimizer for \eqref{e:varopt} is to use \[\vect{\tilde{v}}=\vect{v} + (\beta_v\D{A})^{-1}\D{K}[\iut{\lambda\igrad m}]\] as a search direction (see, e.g.,~\cite{Hart:2009a}). We opt for a (Gauss--)Newton--Krylov method instead, due to its superior rate of convergence (see~\cite{Mang:2015a} for a comparison). Formally, this requires second variations of $\D{L}$. The expression for the action of the reduced space Hessian $\D{H}$ on a vector $\vect{\tilde{v}}$ is given by
\begin{align}
\D{H}[\vect{\tilde{v}}](\vect{v}) \defeq
  \beta_v\D{A}[\vect{\tilde{v}}] + \D{K}[\vect{\tilde{b}}]
= \beta_v\D{A}[\vect{\tilde{v}}] + \D{K}[\iut{\tilde{\lambda} \igrad m + \lambda \igrad \tilde{m}}].
\label{e:hessian-matvec}
\end{align}

The operators $\D{A}$ and $\D{K}$ are as defined above. We, likewise to the reduced gradient $\vect{g}$ in \eqref{e:control-elim}, need to find two space-time fields $\tilde{m}$ and $\tilde{\lambda}$. We can find the \emph{incremental state variable} $\tilde{m}$ by solving
\begin{subequations}
\begin{align}
\p_t \tilde{m} + \vect{v}\cdot\igrad\tilde{m}
+ \vect{\tilde{v}}\cdot\igrad m & = 0
&&\text{in}\;\; \Omega \times (0,1],
\label{e:inc-state-pde-elim}
\\
\tilde{m} &= 0
&&\text{in}\;\;\Omega\times\{0\},
\label{e:inc-state-ic-elim}
\end{align}
\end{subequations}

\noindent with periodic boundary conditions on $\p\Omega$, \emph{forward in time}. Once we have found $\tilde{m}$ we can compute the \emph{incremental adjoint variable} $\tilde{\lambda}$ by solving
\begin{subequations}
\begin{align}
-\p_t \tilde{\lambda} - \idiv(\tilde{\lambda}\vect{v}
+ \lambda\vect{\tilde{v}}) & = 0
&&\text{in}\;\; \Omega \times [0,1),
\label{e:inc-adj-pde-elim}
\\
\tilde{\lambda}  &= -\tilde{m}
&&\text{in}\;\;\Omega\times\{1\},
\label{e:inc-adj-fc-elim}
\end{align}
\end{subequations}

\noindent with periodic boundary conditions on $\p\Omega$, \emph{backward in time}. Thus, each time we apply the Hessian to a vector we have to solve two PDEs---\eqref{e:inc-state-pde-elim} and~\eqref{e:inc-adj-pde-elim}.

\subsection{Discretization}
\label{s:discretization}

We subdivide the time interval $[0,1]$ into $n_t\in\ns{N}$ uniform steps $t^j$, $j=0,\ldots,n_t$, of size $h_t = 1/n_t$. We discretize $\Omega\defeq(-\pi,\pi)^d$ via a regular grid with cell size $\vect{h}_x = (h^1_x,\ldots,h^d_x)^\T\in\ns{R}^d_{>0}$, $\vect{h}_x=2\pi\oslash\vect{n}_x$, $\vect{n}_x=(n^1_x,\ldots,n^2_x)^\T\in\ns{N}^d$; we use a pseudospectral discretization with a Fourier basis. We discretize the integral operators based on a midpoint rule. We use cubic splines as a basis function for our interpolation model.

\subsection{Numerical Time Integration}
\label{s:timeintegration}

An efficient, accurate, and stable time integration of the hyperbolic PDEs that appear in our optimality system is critical for our solver to be effective. Each evaluation of the objective functional $\D{J}$ in~\eqref{e:varopt:objective} requires the solution of~\eqref{e:state-pde-elim} (forward in time). The evaluation of the reduced gradient $\vect{g}$ in~\eqref{e:control-elim} requires an additional solution of~\eqref{e:adj-pde-elim} (backward in time). Applying the reduced space Hessian $\D{H}$ (Hessian matvec) in~\eqref{e:hessian-matvec} necessitates the solution of~\eqref{e:inc-state-pde-elim} (forward in time) and~\eqref{e:inc-adj-pde-elim} (backward in time).

\subsubsection{Second order Runge-Kutta  Schemes}
\label{s:rk2schemes}

In our original work~\cite{Mang:2016a,Mang:2015a} we solved the transport equations based on an RK2 scheme (in particular, Heun's method). This method---in combination with a pseudospectral discretization in space---offers high accuracy solutions, minimal numerical diffusion, and spectral convergence for smooth problems at the cost of having to use a rather small time step due to its conditional stability; the time step size $h_t$ has to be chosen according to considerations of stability rather than accuracy. This scheme can become unstable, even if we adhere to the conditional stability (see \secref{s:experiments:pdesolver} for examples). One strategy to stabilize our solver is to rewrite the transport equations in antisymmetric form~\cite{Fornberg:1975a,Kreiss:1972a}. Here we extend this stable scheme to the adjoint problem and the Hessian operator. We do so by deriving the antisymmetric form of the forward operator and then formally computing its variations. It is relatively straightforward but we have not seen this in the literature related to inverse transport problems. We present the associated PDE operators in \secref{s:rk2a-derivation}. We refer to this solver as RK2A scheme. It is evident that the discretization in antisymmetric from requires more work (see~\secref{s:rk2a-derivation}). We provide estimates in terms of the number of FFTs we have to perform in~\tabref{t:computational-complexity-theoretical} in~\secref{s:computational-complexity}.

\subsubsection{Semi-Lagrangian Formulation}
\label{s:slscheme}

Next, we describe our semi-Lagrangian formulation. To be able to apply the semi-Lagrangian method to the transport equations appearing in our optimality systems, we have to reformulate them. Using the identity $\idiv u\vect{v} = u \idiv \vect{v} + \igrad u\cdot\vect{v}$ for some arbitrary scalar function $u:\bar{\Omega}\rightarrow\ns{R}$, we obtain
\begin{subequations}
\label{e:hyperbolic-pdes-slm}
\begin{align}
\p_t m + \vect{v}\cdot\igrad m & = 0
&&\text{in}\;\;\Omega\times(0,1],
\label{e:state-slm}
\\
- \p_t\lambda
- \vect{v}\cdot\igrad\lambda
- \lambda\idiv\vect{v} & = 0
&&\text{in}\;\;\Omega\times[0,1),
\label{e:adj-slm}
\\
\p_t\tilde{m}
+ \vect{v}\cdot\igrad\tilde{m}
+ \vect{\tilde{v}}\cdot\igrad m & = 0
&&\text{in}\;\;\Omega\times(0,1],
\label{e:inc-state-slm}
\\
- \p_t\tilde{\lambda}
- \vect{v}\cdot\igrad\tilde{\lambda}
- \tilde{\lambda}\idiv\vect{v}
- \idiv\lambda\vect{\tilde{v}} & = 0
&&\text{in}\;\;\Omega\times[0,1).
\label{e:inc-adj-slm}
\end{align}
\end{subequations}

\noindent These equations are all of the general form $\ld u = \p_t u + \vect{v}\cdot\igrad u  = f(u,\vect{v})$, where $u:\bar{\Omega}\times[0,1]\rightarrow\ns{R}$ is some arbitrary scalar function and $\ld \defeq \p_t  + \vect{v}\cdot \igrad$. If the Lagrangian derivative vanishes, i.e., $\ld u = 0$, $u$ is constant along the characteristics $\vect{X} : [\tau_0, \tau_1] \rightarrow \ns{R}^d$ of the flow, where $[\tau_0,\tau_1]\subseteq[0,1]$. We can compute $\vect{X}$ by solving the ODE
\begin{subequations}
\label{e:characteristic}
\begin{align}
\ld \vect{X}(t)
&
=\vect{v}(\vect{X}(t))
&& \text{in}\;\;(\tau_0,\tau_1],
\\
\vect{X}(t)
&
=\vect{x}
&& \text{at}\;\;\{\tau_0\}.
\end{align}
\end{subequations}

\noindent The solution of~\eqref{e:characteristic} requires the knowledge of the velocity field $\vect{v}$ at points that do not coincide with the computational grid; we have to interpolate $\vect{v}$ in space.\footnote{Notice that the scheme becomes more complicated if $\vect{v}$ is non-stationary; we have to interpolate in time and space.}

The idea of pure Lagrangian schemes is to solve $\ld u = f(u,\vect{v})$ along the characteristic lines~\eqref{e:characteristic}. The key advantage of these methods is that they are essentially unconditionally stable~\cite{Staniforth:1991a}; i.e., the time step $h_t$ may be chosen according to accuracy considerations rather than stability considerations.\footnote{For rapidly varying velocity fields, instabilities may still occur.} On the downside the solution will no longer live on a regular grid; the grid changes over time and eventually might become highly irregular. Semi-Lagrangian methods can be viewed as a hybrid between Lagrangian and Eulerian methods; they combine the best from both worlds---they operate on a regular grid and are unconditionally stable.

The semi-Lagrangian scheme involves two steps: For each time step $t^j$ we have to compute the departure point $\vect{X}_D \defeq \vect{X}(t = t^{j-1})$ of a fluid parcel by solving the characteristic equation~\eqref{e:characteristic} backward in time, with initial condition $\vect{X}(t=t^j)=\vect{x}$.\footnote{The direction of time integration depends on the transport equation. For simplicity, we will limit the description of the semi-Lagrangian method to transport equations that are solved forward in time. Notice that \eqref{e:hyperbolic-pdes-slm} also contains equations that have to be solved backward in time.} We revert to a uniform grid by interpolation. The second step is to compute the transported quantity along the characteristic $\vect{X}$. The accuracy of the semi-Lagrangian method is sensitive to the time integrator for solving~\eqref{e:characteristic} as well as the interpolation scheme used to evaluate the departure points $\vect{X}_D$. We discuss the individual building blocks of our solver next.

\paragraph{Tracing the Characteristic}

\begin{figure}
\centering
\includegraphics[width=0.5\textwidth]
{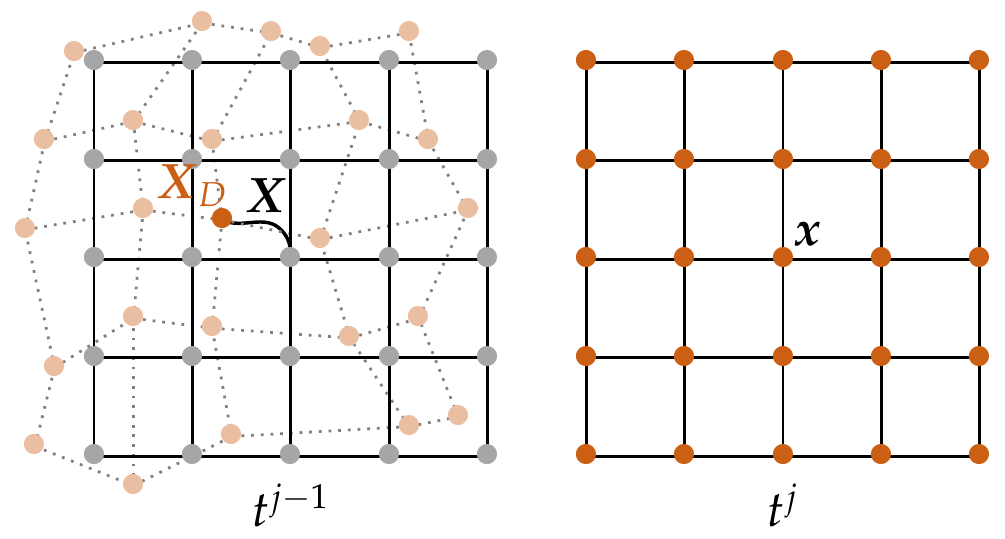}
\caption{Tracing the characteristic $\vect{X}$ in a semi-Lagrangian scheme. We start with a regular grid $\Omega^h$ (dark orange points on the right) consisting of coordinates $\vect{x}$ at time point $t^j$. We assume we have already computed the intermediate solution $u^h$ of a given transport equation at time point $t^{j-1}$; we know the input data on the regular grid at time point $t^{j-1}$ (the regular grid nodes are illustrated in light gray). In a first step, we trace back the characteristic $\vect{X}$ by solving~\eqref{e:characteristic} backward in time subject to the initial condition $\vect{X}(t=t^j) = \vect{x}$. Once we have found the characteristic (black line in the figure on the left) we can---in a second step---assign the value of $u^h$ at $t^j$ given at the departure point $\vect{X}_D$ (dark orange point on the left) to $\vect{x}$ at $t^j$ based on some interpolation model. We illustrate the grid of departure points in light orange and the original grid in gray (left figure).
\label{f:tracing-the-characteristic}}
\end{figure}

For each time step $t^j$ of the integration of a given transport equation we have to trace the characteristic $\vect{X}$ backward in time in an interval $[t^{j-1},t^j]\subset[0,1]$. Since we invert for a stationary velocity field $\vect{v}$ we have to trace $\vect{X}$ (i.e., compute the departure points $\vect{X}_D$) only once in every Newton iterations used for all time steps.\footnote{In total, we actually need to compute two characteristics, one for the forward (state or primal) equations and one for the backward (adjoint or dual) equations.} We use an explicit RK2 scheme (Heun's method) to do so~\cite{Staniforth:1991a}. We illustrate the computation of the characteristic in \figref{f:tracing-the-characteristic}. Each evaluation of the right hand side of~\eqref{e:characteristic} requires interpolation.

\paragraph{Interpolation}

We use a cubic spline interpolation model to evaluate the transported quantities along the characteristic $\vect{X}$. We pad the imaging data to account for the periodic boundary conditions. The size of the padding zone is computed at every iteration based on the maximal displacement between the original grid nodes and the departure points $\vect{X}_D$; we also account for the support of the basis functions of the interpolation model.\footnote{We have tested a more accurate implementation that applies the interpolation step on a grid of half the cell size $\vect{h}_x$ to minimize the interpolation error as well as the numerical diffusion. We prolong and restrict the data in the Fourier domain. The gain in numerical accuracy (one to two digits) did not justify the significant increase in CPU time.}

\paragraph{Transport}

To transport the quantity of interest we have to solve equations of the form
\begin{align}
\label{e:transport-along-characteristic}
\ld u(\vect{X}(t),t)
&
= f(u(\vect{X}(t)),t),\vect{v}(\vect{X}(t)))
&&\text{in}\;\; [t^{j-1}, t^j]
\end{align}

\noindent along the characteristic $\vect{X}$. We use an explicit RK2 scheme to numerically solve~\eqref{e:transport-along-characteristic}. Since $u$ will be needed along the characteristic $\vect{X}$ we have to interpolate $u$ at the computed departure points $\vect{X}_D$.

\subsection{Numerical Optimization}
\label{s:newtonkrylovmethod}

We use a globalized, inexact, matrix-free (Gauss--)Newton--Krylov method for numerical optimization. Our solver has been described and tested in~\cite{Mang:2015a}. In what follows, we will briefly revisit this solver and from thereon design a nested, two-level preconditioner for the reduced space optimality conditions.

\subsubsection{Newton--Krylov Solver}

The Newton step for updating $\vect{v}^h_k\in\ns{R}^n$, $n=d\prod_{i=1}^d n_x^i$, is in general format given by
\begin{equation}
\label{e:newtonsteprs}
\D{H}^h\vect{\tilde{v}}^h_k = -\vect{g}^h_k,
\qquad
\vect{v}^h_{k+1} = \vect{v}^h_k
+ \alpha_k\tilde{\vect{v}}^h_k,
\end{equation}

\noindent where $\D{H}^h\in\ns{R}^{n, n}$, is the reduced space Hessian operator, $\vect{\tilde{v}}^h_k\in\ns{R}^n$ the search direction, and $\vect{g}^h_k\in\ns{R}^n$ the reduced gradient.\footnote{The Hessian matvec is given by \eqref{e:hessian-matvec}; the expression for the reduced gradient is given in \eqref{e:control-elim}.} We globalize our iterative scheme on the basis of a backtracking line search subject to the Armijo--Goldstein condition with step size $\alpha_k>0$ at iteration $k\in\ns{N}$ (see, e.g.,~\cite[page~37]{Nocedal:2006a}). We keep iterating until the relative change of the gradient $\|\vect{g}^h_k\|_{\text{rel}}\defeq\|\vect{g}^h_k\|_{\infty}/\|\vect{g}^h_0\|_{\infty}$ is smaller or equal to $\num{1E-2}$ or $\|\vect{g}^h_k\|_{\infty} \leq \num{1E-5}$ (other stopping conditions can be used; see e.g.~\cite[pages~305\,ff.]{Gill:1981a}). We refer to the steps necessary for updating $\vect{v}^h_k$ as \emph{outer iterations} and to the steps necessary for ``inverting the reduced Hessian'' in~\eqref{e:newtonsteprs} (i.e., the steps necessary to solve for the search direction $\vect{\tilde{v}}^h_k$) as \emph{inner iterations} (see \cite{Mang:2016a,Mang:2015a} for more details).

We use a Krylov iterative solver to compute $\vect{\tilde{v}}^h_k$. To evaluate the reduced gradient $\vect{g}^h_k\in\ns{R}^n$ on the right hand side of~\eqref{e:newtonsteprs} (see~\eqref{e:control-elim}) we have to solve~\eqref{e:state-pde-elim} \emph{forward in time} and~\eqref{e:adj-pde-elim} \emph{backward in time} for a given iterate $\vect{v}^h_k$.\footnote{In general we associate the cost for   solving~\eqref{e:state-pde-elim} to the evaluation of the objective in~\eqref{e:varopt:objective}, i.e., to the line search.} Once we have found the gradient, we can solve~\eqref{e:newtonsteprs}. The reduced space Hessian in~\eqref{e:newtonsteprs} is a large, dense, ill-conditioned operator. Solving this system is a significant challenge; we use a PCG method~\cite{Hestenes:1952a}. Indefiniteness of $\D{H}^h$ can be avoided by using a GN approximation to the true Hessian\footnote{This corresponds to dropping all terms with $\lambda$ in~\eqref{e:hessian-matvec} and~\eqref{e:inc-adj-pde-elim} (see~\cite{Mang:2016a,Mang:2015a}).} or by terminating the PCG solve in case negative curvature occurs. By using a GN approximation, we sacrifice speed of convergence; quadratic convergence drops to superlinear convergence; we locally recover quadratic convergence as $\lambda^h$ tends to zero.

An important property of Krylov subspace methods is that we do not have to store or form the reduced space Hessian $\D{H}^h$; we merely need an expression for the action of $\D{H}^h$ on a vector; this is exactly what~\eqref{e:hessian-matvec} provides. Each application of $\D{H}^h$ (i.e., each PCG iteration) requires the solution of~\eqref{e:inc-state-pde-elim} and~\eqref{e:inc-adj-pde-elim}. This results in high computational costs. We use inexact solves (see~\cite[pages\,165ff.]{Nocedal:2006a} and references therein) to reduce these costs. Another key ingredient to keep the number of PCG iterations small is an effective preconditioner. This is what we discuss next.

\subsubsection{Preconditioner}
\label{s:precond}

The design of an optimal preconditioner for KKT systems arising in large-scale inverse problems is an active area of research~\cite{Benzi:2011a,Biros:2008a,Biros:2005a,Biros:2005b,Haber:2001a}.\footnote{We study the spectral properties of $\D{H}^h$ for the compressible and incompressible case in~\cite{Mang:2015a}.} Standard techniques, like incomplete factorizations, are not applicable as they require the assembling of $\D{H}^h$. We provide two matrix-free strategies below.

Given the reduced gradient $\vect{g}$, the Newton step in the reduced space is given by
\begin{equation}
\label{e:reduced-space-kkt-system}
\D{H}[\vect{\tilde{v}}](\vect{v})
=
\beta_v\D{A}[\vect{\tilde{v}}] +
\D{K}
\left[
\int_0^1 \tilde{\lambda}\igrad m + \lambda\igrad \tilde{m} \d{t}
\right]
= \beta_v\D{A}[\vect{\tilde{v}}] + \D{Q}[\vect{\tilde{v}}]
= -\vect{g}.
\end{equation}

\noindent We have introduced the operator $\D{Q}[\vect{\tilde{v}}] \defeq\D{Q}[\tilde{\lambda}, m, \lambda, \tilde{m}](\vect{v},\vect{\tilde{v}})$ in~\eqref{e:reduced-space-kkt-system} for notational convenience and to better illustrate its dependence on $\vect{\tilde{v}}$; the incremental state and adjoint variables, $\tilde{m}$ and $\tilde{\lambda}$, are functions of $\vect{\tilde{v}}$ through~\eqref{e:inc-state-pde-elim} and~\eqref{e:inc-adj-pde-elim}, respectively.

We use a left preconditioner $\mat{P}^{-1}$; our solver will see the system $\mat{P}^{-1}\D{H}^h\vect{\tilde{v}}^h_k = -\mat{P}^{-1}\vect{g}^h_k$. Ideally the preconditioned matrix will have a much better spectral condition number and/or eigenvalues that are clustered around one. An ideal preconditioner is one that has vanishing costs for its construction and application and at the same time represents an excellent approximation to the Hessian operator $\D{H}^h$ so that $\mat{P}^{-1}\D{H}^h\approx \mat{I}_n$~\cite{Benzi:2002a}. These are in general competing goals. Since we use a PCG method to iteratively solve~\eqref{e:reduced-space-kkt-system}, we only require the action of $\mat{P}^{-1}$ on a vector.

\paragraph{Regularization Preconditioner}

In our original work~\cite{Mang:2016a,Mang:2015a}, we use a preconditioner that is based on the exact, spectral inverse of the regularization operator $\D{A}^h$, i.e.,
\begin{equation}
\label{e:reg-precond}
\mat{P}_{\text{REG}}
= \beta_v\D{A}^h
= \beta_v\mat{W}\mat{\Gamma}\mat{W}^{-1},
  \quad \mat{P}_{\text{REG}}\in\ns{R}^{n, n},
\end{equation}

\noindent where $\mat{W}^{-1} = \mat{I}_d\otimes\hat{\mat{W}}\in\ns{C}^{n,n}$, $\mat{I}_d = \diag(1,\ldots,1)\in\ns{R}^{d,d}$, $\hat{\mat{W}}$ is a DFT matrix and $\mat{\Gamma} = \mat{I}_d\otimes\hat{\mat{\Gamma}}\in\ns{R}^{n,n}$ are the spectral weights for the Laplacian, biharmonic, or triharmonic differential operators in~\eqref{e:diffopregv}. The operator $\D{A}^h$ has a non-trivial kernel; to be able to invert this operator analytically we replace the zero entries in $\mat{\Gamma}$ by one. If we apply $\mat{P}_{\text{REG}}^{-1}$ to the reduced Hessian in~\eqref{e:reduced-space-kkt-system} the system we are effectively solving is a low-rank, compact perturbation of the identity:
\begin{equation}
\label{e:reg-precond-hessian}
\vect{\tilde{v}}^h + (\beta_v\D{A}^h)^{-1}
\D{Q}^h[\vect{\tilde{v}}^h]
= (\mat{I}_n+(\beta_v\D{A}^h)^{-1}\D{Q}^h)
\vect{\tilde{v}}^h.
\end{equation}

\noindent Notice that the operator $\mat{P}_{\text{REG}}$ acts as a smoother on $\D{Q}^h$. Applying and inverting this preconditioner has vanishing computational costs (due to our pseudospectral discretization). This preconditioner becomes ineffective for small regularization parameters $\beta_v$ and a high inversion accuracy (i.e., small tolerances for the relative reduction of the reduced gradient; see, e.g.,~\cite{Mang:2016a}).

\paragraph{Nested Preconditioner}

We use a coarse grid correction by an inexact solve to provide an improved preconditioner. This corresponds to a simplified two-level multigrid v-cycle, where the smoother is the inverse of the regularization operator and the coarse grid solve is inexact. We introduce spectral restriction and prolongation operators to change from the fine to the coarse grid and vice versa. The action of the preconditioner, i.e., of the action of the reduced space Hessian in~\eqref{e:reduced-space-kkt-system}, is computed on the coarse grid. This preconditioner only operates on the low frequency modes due to the restriction to a coarser grid. In our implementation, we treat the high and the low frequency components separately; we apply the nested preconditioner to the low frequency modes and leave the high frequency modes untouched. We separate the frequency components by applying an ideal low- and high-pass filter to the vector the preconditioner is applied to. As we will see below, we will actually treat the high frequency components with a smoother that is based on the inverse of our regularization operator, i.e., the Hessian to be preconditioned does not correspond to the reduced space Hessian in~\eqref{e:reduced-space-kkt-system} but the preconditioned Hessian in~\eqref{e:reg-precond-hessian}. We refer to this preconditioner as $\mat{P}_{\text{2L}}$.

The effectiveness of this scheme is dictated by the computational costs associated with the inversion of the (coarse grid) Hessian operator $\D{H}^{2h}\in\ns{R}^{n/2,n/2}$. One strategy for applying this preconditioner is to compute the action of the inverse of $\D{H}^{2h}$ using a nested PCG method. From the theory of Krylov subspace methods we know that we have to solve for the action of this inverse with a tolerance that is smaller than the one we use to solve~\eqref{e:reduced-space-kkt-system} (exact solve; we refer to this approach as PCG($\epsilon$), where $\epsilon\in(0,1)$ is the scaling for the tolerance used to solve~\eqref{e:reduced-space-kkt-system}) for the outer PCG method to not break down. This increased accuracy may lead to impractical computational costs, especially since each application of $\D{H}^{2h}$ requires two PDE solves, and we expect this preconditioner to have a very similar conditioning as $\D{H}^h$. Another strategy is to solve the system inexactly. This requires the use of flexible Krylov subspace methods (for the Hessian operator) or a Chebyshev semi-iterative method (CHEB; for the preconditioner) (see e.g.~\cite[pages\,179ff.]{Axelsson:1996a}) with a fixed number of iterations (we refer to this strategy as CHEB($k$), where $k$ is the number of iterations). This makes the work spent on inverting the preconditioner constant but the inexactness might lead to a less effective preconditioner. Another bottleneck is the fact that the CHEB method requires estimates of the spectral properties of of the operator we try to invert; estimating the eigenvalues is expensive and can lead to excessive computational costs. We provide implementation details next, some of which are intended to speed up the formation and application of our nested preconditioner.

\begin{itemize}
\item \textbf{Spectral Preconditioning of $\D{H}^h$}: Since the application of the inverse of the regularization operator $\D{A}^h$ comes at almost no cost, we decided to use the spectrally preconditioned Hessian operator in~\eqref{e:reg-precond-hessian} within our two-level scheme, with a small technical modification. The left preconditioned Hessian in~\eqref{e:reg-precond-hessian} is not symmetric. We can either opt for Krylov methods that do not require the operator we try to invert to be symmetric, or we employ a spectral split preconditioner. We opt for the latter approach to be able to use a PCG method, attributed to its efficiency. The split preconditioned system is given by \[(\mat{I}_n + (\beta_v\D{A})^{-1/2}\D{Q}^h (\beta_v\D{A})^{-1/2})\vect{s} = -(\beta_v\D{A})^{-1/2}\vect{g},\] where $\vect{s}\defeq (\beta_v\D{A})^{1/2}\vect{\tilde{v}}$. Notice, that the inverse of the regularization operator can be viewed as a smoother, which establishes a connection of our scheme to more sophisticated multigrid strategies~\cite{Adavani:2008a,Borzi:2002a}.
\item \textbf{Eigenvalue Estimates for the CHEB Method}: The computational costs for estimating eigenvalues of $\mat{P}_{\text{2L}}$ are significant. Our assumption is that we have to estimate the extremal eigenvalues only once for the registration for a given set of images (we will experimentally verify this assumption; see \secref{s:experiments:eigenvalest}); if we change the regularization parameter we simply have to scale the estimated eigenvalues. Notice that we can efficiently estimate the eigenvalues for a zero velocity field since a lot of the terms drop in the optimality systems. We compute an estimate for the largest eigenvalue $e_{\max}$ based on an implicitly restarted Lanczos algorithm. We approximate the smallest eigenvalue analytically under the assumption that $\D{Q}^h$ is a low-rank operator of order $\bigO(1)$; $e_{\min} \approx \min(\mat{I}_n + (\beta_v\mat{\Gamma})^{-1})$.
\item \textbf{Hyperbolic PDE Solves}: Each matvec with $\D{H}^{2h}$ requires the solution of~\eqref{e:inc-state-pde-elim} and~\eqref{e:inc-adj-pde-elim} on the coarse grid. We exclusively consider the semi-Lagrangian formulation to speed up the computation. In general we assume that we do not need high accuracy solutions for our preconditioner. This might even be true for the PDE solves within each Hessian matvec.\footnote{We have also tested an approximation of the forcing term by dropping all second order terms of the RK2 scheme for numerically integrating~\eqref{e:inc-state-pde-elim} and~\eqref{e:inc-adj-pde-elim}. Since we have observed instabilities in the RK2 schemes and due to the effectiveness of the semi-Lagrangian method (see \secref{s:experiments:pdesolver}) we do not report results for this preconditioner.}
\item \textbf{Restriction}/\textbf{Prolongation}: We use spectral restriction and prolongation operators. We do not apply an additional smoothing step after or before we restrict the data to the coarser grid. We actually observed that applying an additional Gaussian smoothing with a standard deviation of $2\vect{h}_x$ (i.e., one grid point on the coarser grid) significantly deteriorates the performance of our preconditioner for small grid sizes (e.g., $64\times64$). A more detailed study on how the choices for the restriction and prolongation operators affect the performance of our solver with respect to changes in the regularity of the underlying objects remains for future work.
\item \textbf{Filters}: We use simple cut-off filters before applying the restriction and prolongation operators, with a cut-off frequency of half the frequency that can be represented on the finer grid.
\end{itemize}

\subsection{Implementation Details and Parameter Settings}

Here, we briefly summarize some of the implementation details and parameter choices.
\begin{itemize}
\item \textbf{Image Data}: Our solver can not handle images with discontinuities. We ensure that the images are adequately smooth by applying a Gaussian smoothing kernel with an empirically selected standard deviation of one grid point in each spatial direction. We normalize the intensities of the images to $[0,1]$ prior to registration.
\item \textbf{PDE Solves}: We use a CFL number of 0.2 for the explicit RK2 schemes; we observed instabilities for some of the test cases for a CFL number of 0.5. The semi-Lagrangian method is unconditionally stable; we test different CFL numbers.
\item \textbf{Restriction}/\textbf{Prolongation}: We use spectral prolongation and restriction operators within our preconditioner (more implementation details for our preconditioner can be found in the former section). We do not perform any other grid, scale, or parameter continuation to speed up our computations.
\item \textbf{Interpolation}: We consider a $C^2$-continuous cubic spline interpolation model. We extend our data periodically to account for the boundary conditions.
\item \textbf{Regularization}: Since we study the behavior of our solver as a function of the regularization parameters, we will set their value empirically. For practical applications, we have designed a strategy that allows us to probe for an ideal regularization parameter; we perform a parameter continuation that is based on a binary search and considers bounds on the determinant of the deformation gradient as a criterion; see~\cite{Mang:2015a,Mang:2016a}.
\item \textbf{Globalization}: We use a backtracking line search subject to the Armijo--Goldstein condition to globalize our Newton--Krylov scheme (see, e.g.,~\cite[page~37]{Nocedal:2006a}).
\item \textbf{Stopping Criteria}: We terminate the inversion if the relative change of the gradient is smaller or equal to $\num{1E-2}$ or $\|\vect{g}^h_k\|_{\infty} \leq \num{1E-5}$  (other stopping conditions can be used; see, e.g.,~\cite[pages~305\,ff.]{Gill:1981a}).
\item \textbf{Hessian}: We use a GN approximation to the reduced space Hessian $\D{H}^h$ to avoid indefiniteness. This corresponds to dropping all expressions with $\lambda$ in~\eqref{e:hessian-matvec} and~\eqref{e:inc-adj-pde-elim} (see~\cite{Mang:2015a} for more details); we recover quadratic convergence for $\lambda\rightarrow0$.
\item \textbf{KKT solve}: If not noted otherwise, we will solve the reduced space KKT system in~\eqref{e:newtonsteprs} inexactly, with a forcing sequence that assumes quadratic convergence (see~\cite[pages\,165ff.]{Nocedal:2006a} and references therein); we use a PCG method to iteratively solve~\eqref{e:newtonsteprs}.
\item \textbf{PC solve}: We compute the action of the inverse of the 2-level preconditioner either exactly using a nested PCG method or inexactly based on a nested CHEB method with a fixed number of iterations.
\end{itemize}

\begin{landscape}
\resetrunid
\begin{table}
\caption{Self-convergence for the RK2, the RK2A, and the SL method for the numerical integration of the state (see~\eqref{e:state-pde-elim}; results reported in top block) and the adjoint (see~\eqref{e:adj-pde-elim}; results reported in bottom block) equation. We report the relative $\ell^2$-error $\|\delta u^h\|_{\text{rel}}$ between solutions for the state ($u^h = m_1^h$) and the adjoint ($u^h=\lambda_0^h$) equation computed on a grid of size $\vect{n}_x$ and a grid of size $2\vect{n}_x$. We use a CFL number of 0.2 for the RK2 and the RK2A method, and a CFL number of 0.2, 1, and 5, for the SL method; we provide the associated number of time points $n_t$. We report errors for different grid sizes and test problems (top block: SMOOTH A; bottom block: SMOOTH B; see \figref{f:synthetic-data}); $\ast\ast\ast$ indicates that the solver became unstable (not due to a violation of the CFL condition; see text for details). We also report the time to solution (in seconds).
\label{t:selfconvergence-hyperbolic-pdesolvers}}
\centering
\begin{tiny}
\begin{tabular}{rrrrrLlrrLlrrLlrrLlrrLl}
\toprule
& & & \mcol{4}{RK2(0.2)} & \mcol{4}{RK2A(0.2)} & \mcol{4}{SL(0.2)} & \mcol{4}{SL(1)} & \mcol{4}{SL(5)}
\\\midrule
& & $n_x^i$ & run & $n_t$ & $\|\delta u^h\|_{\text{rel}}$ & time & run & $n_t$ & $\|\delta u^h\|_{\text{rel}}$ & time & run & $n_t$ & $\|\delta u^h\|_{\text{rel}}$ & time & run & $n_t$ & $\|\delta u^h\|_{\text{rel}}$ & time & run & $n_t$ & $\|\delta u^h\|_{\text{rel}}$ & time
\\\midrule
  \mrowrot{8}{STATE EQ}
& A & 64
& \runid &  26 & \num{1.246470e-05} & \num{3.249150e-01}
& \runid &  26 & \num{1.246470e-05} & \num{6.642610e-01}
& \runid &  26 & \num{3.790105e-06} & \num{8.526410e-01}
& \runid &   6 & \num{4.931513e-05} & \num{2.696060e-01}
& \runid &   2 & \num{3.185243e-04} & \num{5.875400e-02}
\\
& & 128
& \runid &  51 & \num{3.279612e-06} & \num{1.629233e+00}
& \runid &  51 & \num{3.279612e-06} & \num{3.201515e+00}
& \runid &  51 & \num{8.628094e-07} & \num{3.849134e+00}
& \runid &  11 & \num{1.518920e-05} & \num{1.316970e+00}
& \runid &   3 & \num{1.685669e-04} & \num{2.408120e-01}
\\
& & 256
& \runid & 102 & \num{8.194718e-07} & \num{1.069440e+01}
& \runid & 102 & \num{8.194718e-07} & \num{2.040621e+01}
& \runid & 102 & \num{2.000567e-07} & \num{3.038315e+01}
& \runid &  21 & \num{4.255715e-06} & \num{6.809877e+00}
& \runid &   5 & \num{6.752966e-05} & \num{1.564274e+00}
\\
& & 512
& \runid & 204 & \num{2.048135e-07} & \num{8.440371e+01}
& \runid & 204 & \num{2.048135e-07} & \num{1.609808e+02}
& \runid & 204 & \num{4.807773e-08} & \num{2.525976e+02}
& \runid &  41 & \num{1.138702e-06} & \num{5.087391e+01}
& \runid &   9 & \num{2.215575e-05} & \num{1.141639e+01}
\\\cmidrule{2-23}
& B & 64
& \runid & 102 & \num{2.362605e-01} & \num{9.204050e-01}
& \runid & 102 & \num{7.669804e-02} & \num{2.382031e+00}
& \runid & 102 & \num{9.934129e-02} & \num{3.331444e+00}
& \runid &  21 & \num{8.635175e-02} & \num{6.780110e-01}
& \runid &   5 & \num{7.351601e-02} & \num{2.235910e-01}
\\
& & 128
& \runid & 204 & \num{2.812852e-02} & \num{6.440490e+00}
& \runid & 204 & \num{9.832740e-03} & \num{1.262654e+01}
& \runid & 204 & \num{1.057699e-02} & \num{1.772764e+01}
& \runid &  41 & \num{9.670545e-03} & \num{4.598787e+00}
& \runid &   9 & \num{9.810264e-03} & \num{1.153501e+00}
\\
& & 256
& \runid & 408 & $\ast\ast\ast$     & \num{4.228000e+01}
& \runid & 408 & \num{1.199602e-04} & \num{7.855845e+01}
& \runid & 408 & \num{4.027374e-04} & \num{1.161691e+02}
& \runid &  82 & \num{3.671445e-04} & \num{2.513303e+01}
& \runid &  17 & \num{1.192550e-03} & \num{6.938980e+00}
\\
& & 512
& \runid & 815 & $\ast\ast\ast$     & \num{3.431830e+02}
& \runid & 815 & \num{3.875719e-06} & \num{6.522741e+02}
& \runid & 815 & \num{1.875843e-05} & \num{1.060886e+03}
& \runid & 163 & \num{2.743683e-05} & \num{1.958110e+02}
& \runid &  33 & \num{3.095846e-04} & \num{4.225729e+01}
\\\midrule
  \mrowrot{8}{ADJOINT EQ}
& A & 64
& \runid &  26 & \num{1.284931e-04} & \num{2.136060e-01}
& \runid &  26 & \num{1.284931e-04} & \num{6.598990e-01}
& \runid &  26 & \num{8.334623e-05} & \num{7.330560e-01}
& \runid &   6 & \num{1.389942e-03} & \num{3.053540e-01}
& \runid &   2 & \num{8.360937e-03} & \num{1.637680e-01}
\\
& & 128
& \runid &  51 & \num{3.398516e-05} & \num{1.355887e+00}
& \runid &  51 & \num{3.398516e-05} & \num{3.416255e+00}
& \runid &  51 & \num{2.249425e-05} & \num{6.041776e+00}
& \runid &  11 & \num{4.469800e-04} & \num{1.302645e+00}
& \runid &   3 & \num{4.606136e-03} & \num{5.304680e-01}
\\
& & 256
& \runid & 102 & \num{8.515454e-06} & \num{9.287153e+00}
& \runid & 102 & \num{8.515454e-06} & \num{2.362273e+01}
& \runid & 102 & \num{5.662382e-06} & \num{3.575260e+01}
& \runid &  21 & \num{1.280244e-04} & \num{8.235509e+00}
& \runid &   5 & \num{1.928118e-03} & \num{2.082276e+00}
\\
& & 512
& \runid & 204 & \num{2.131305e-06} & \num{7.094955e+01}
& \runid & 204 & \num{2.131305e-06} & \num{1.625894e+02}
& \runid & 204 & \num{1.418703e-06} & \num{2.917081e+02}
& \runid &  41 & \num{3.464365e-05} & \num{6.368748e+01}
& \runid &   9 & \num{6.535220e-04} & \num{1.564462e+01}
\\\cmidrule{2-23}
& B & 64
& \runid & 102 & \num{4.043174e-01} & \num{8.209290e-01}
& \runid & 102 & \num{4.187172e-01} & \num{2.234389e+00}
& \runid & 102 & \num{6.066673e-01} & \num{3.274474e+00}
& \runid &  21 & \num{5.299555e-01} & \num{7.785470e-01}
& \runid &   5 & \num{4.342303e-01} & \num{2.824540e-01}
\\
& & 128
& \runid & 204 & \num{5.720269e-02} & \num{6.306302e+00}
& \runid & 204 & \num{5.940875e-02} & \num{1.137439e+01}
& \runid & 204 & \num{6.562716e-02} & \num{1.796776e+01}
& \runid &  41 & \num{6.059585e-02} & \num{3.391594e+00}
& \runid &   9 & \num{5.692348e-02} & \num{1.118745e+00}
\\
& & 256
& \runid & 408 & $\ast\ast\ast$     & \num{4.203301e+01}
& \runid & 408 & \num{7.311170e-04} & \num{7.777944e+01}
& \runid & 408 & \num{2.368267e-03} & \num{1.463573e+02}
& \runid &  82 & \num{2.036229e-03} & \num{2.307441e+01}
& \runid &  17 & \num{3.916051e-03} & \num{4.810801e+00}
\\
& & 512
& \runid & 815 & $\ast\ast\ast$     & \num{2.813443e+02}
& \runid & 815 & \num{1.170932e-05} & \num{6.417905e+02}
& \runid & 815 & \num{1.077794e-04} & \num{1.190277e+03}
& \runid & 163 & \num{1.132407e-04} & \num{2.288566e+02}
& \runid &  33 & \num{1.006256e-03} & \num{5.026806e+01}
\\\bottomrule
\end{tabular}
\end{tiny}
\end{table}

\begin{table}
\caption{Convergence of the SL method to a reference solution computed via the RK2A scheme. We compute the reference solution on a grid of size $\vect{n}_x=(512,512)^\T$ using the RK2A scheme with a CFL number of 0.2. We report results for varying discretization sizes. We report the CFL number $c$, the associated number of time steps $n_t$, and the relative $\ell^2$-error between the solution for the SL scheme and the reference solution computed via the RK2A scheme. We also report errors for the RK2A method as a reference (self convergence). We report results for the state equation (two blocks on the left; see \eqref{e:state-pde-elim}) and the adjoint equation (two blocks on the right; see \eqref{e:adj-pde-elim}). We consider the test problems SMOOTH A and SMOOTH B in \figref{f:synthetic-data} for the velocity field and to set up the initial and terminal conditions, respectively.
\label{t:convergence-SL-to-RK2A-transporteqs}}
\centering
\resetrunid
\begin{tiny}
\begin{tabular}{rrrrRrrrrRrrrrRrrrrRrr}
\toprule
  \mcol{12}{STATE EQ}
& \mcol{10}{ADJOINT EQ}
\\\midrule
  \mcol{7}{SMOOTH A}
& \mcol{5}{SMOOTH B}
& \mcol{5}{SMOOTH A}
& \mcol{5}{SMOOTH B}
\\\midrule
  $n_x^i$ & $c$ & $n_t$ & run & SL & run & RK2A & $n_t$ & run & SL & run & RK2A & $n_t$ & run & SL & run & RK2A & $n_t$ & run & SL & run & RK2A
\\\midrule
  64
& 10
&  2 & \runid & \num{5.838125e-04} && ---
&  3 & \runid & \num{8.106376e-02} && ---
&  2 & \runid & \num{1.570851e-02} && ---
&  3 & \runid & \num{4.178084e-01} && ---
\\
& 5
&  2 & \runid & \num{5.838125e-04} && ---
&  5 & \runid & \num{7.575585e-02} && ---
&  2 & \runid & \num{1.570851e-02} && ---
&  5 & \runid & \num{4.408773e-01} && ---
\\
& 2
&  3 & \runid & \num{2.672315e-04} && ---
& 11 & \runid & \num{7.994546e-02} && ---
&  3 & \runid & \num{7.437041e-03} && ---
& 11 & \runid & \num{4.902524e-01} && ---
\\
& 1
&  6 & \runid & \num{7.043158e-05} && ---
& 21 & \runid & \num{8.566296e-02} && ---
&  6 & \runid & \num{2.010841e-03} && ---
& 21 & \runid & \num{5.307796e-01} && ---
\\
& 0.2
&  26 & \runid & \num{5.103671e-06}
      & \runid & \num{1.656299e-05}
& 102 & \runid & \num{9.811464e-02}
      & \runid & \num{7.663254e-02}
&  26 & \runid & \num{1.139734e-04}
      & \runid & \num{1.709814e-04}
& 102 & \runid & \num{6.064323e-01}
      & \runid & \num{4.188275e-01}
\\\midrule
  128
&  10
&   2 & \runid & \num{5.836446e-04} && ---
&   5 & \runid & \num{1.780858e-02} && ---
&   2 & \runid & \num{1.570732e-02} && ---
&   5 & \runid & \num{7.095293e-02} && ---
\\
&   5
&   3 & \runid & \num{2.668500e-04} && ---
&   9 & \runid & \num{1.078409e-02} && ---
&   3 & \runid & \num{7.434422e-03} && ---
&   9 & \runid & \num{5.911109e-02} && ---
\\
&   2
&   6 & \runid & \num{6.902144e-05} && ---
&  21 & \runid & \num{9.545266e-03} && ---
&   6 & \runid & \num{2.008532e-03} && ---
&  21 & \runid & \num{5.897257e-02} && ---
\\
& 1
&  11 & \runid & \num{2.114234e-05} && ---
&  41 & \runid & \num{9.837339e-03} && ---
&  11 & \runid & \num{6.218182e-04} && ---
&  41 & \runid & \num{6.142297e-02} && ---
\\
& 0.2
&  51 & \runid & \num{1.328577e-06}
      & \runid & \num{4.099048e-06}
& 204 & \runid & \num{1.072697e-02}
      & \runid & \num{9.831807e-03}
&  51 & \runid & \num{3.077680e-05}
      & \runid & \num{4.250007e-05}
& 204 & \runid & \num{6.646499e-02}
      & \runid & \num{5.941213e-02}
\\\midrule
  256
& 10
&   3 & \runid & \num{2.668275e-04} && ---
&   9 & \runid & \num{5.059342e-03} && ---
&   3 & \runid & \num{7.434263e-03} && ---
&   9 & \runid & \num{1.546127e-02} && ---
\\
& 5
&   5 & \runid & \num{9.857960e-05} && ---
&  17 & \runid & \num{1.603810e-03} && ---
&   5 & \runid & \num{2.844454e-03} && ---
&  17 & \runid & \num{5.213297e-03} && ---
\\
& 2
&  11 & \runid & \num{2.096355e-05} && ---
&  41 & \runid & \num{5.016281e-04} && ---
&  11 & \runid & \num{6.214298e-04} && ---
&  41 & \runid & \num{2.086868e-03} && ---
\\
&   1
&  21 & \runid & \num{5.958274e-06} && ---
&  82 & \runid & \num{3.930842e-04} && ---
&  21 & \runid & \num{1.749265e-04} && ---
&  82 & \runid & \num{2.122254e-03} && ---
\\
& 0.2
& 102 & \runid & \num{4.917776e-07}
      & \runid & \num{8.194718e-07}
& 408 & \runid & \num{4.213892e-04}
      & \runid & \num{1.199602e-04}
& 102 & \runid & \num{8.614967e-06}
      & \runid & \num{8.515454e-06}
& 408 & \runid & \num{2.467635e-03}
      & \runid & \num{7.311170e-04}
\\\midrule
  512
& 10
&   5 & \runid & \num{9.857666e-05} && ---
&  17 & \runid & \num{1.473067e-03} && ---
&   5 & \runid & \num{2.844433e-03} && ---
&  17 & \runid & \num{4.597292e-03} && ---
\\
& 5
&   9 & \runid & \num{3.108795e-05} && ---
&  33 & \runid & \num{4.150383e-04} && ---
&   9 & \runid & \num{9.182715e-04} && ---
&  33 & \runid & \num{1.346180e-03} && ---
\\
& 2
&  21 & \runid & \num{5.935902e-06} && ---
&  82 & \runid & \num{8.606744e-05} && ---
&  21 & \runid & \num{1.748711e-04} && ---
&  82 & \runid & \num{2.853485e-04} && ---
\\
& 1
&  41 & \runid & \num{1.711848e-06} && ---
& 163 & \runid & \num{3.593039e-05} && ---
&  41 & \runid & \num{4.696497e-05} && ---
& 163 & \runid & \num{1.285623e-04} && ---
\\
& 0.2
& 204 & \runid & \num{3.193784e-07}
      & \runid & 0
& 815 & \runid & \num{2.086078e-05}
      & \runid & 0
& 204 & \runid & \num{3.729688e-06}
      & \runid & 0
& 815 & \runid & \num{1.049517e-04}
      & \runid & 0
\\\bottomrule
\end{tabular}
\end{tiny}
\end{table}
\end{landscape}

\section{Numerical Experiments}
\label{s:experiments}

We report numerical experiments next. The error of our discrete approximation to the control problem depends on the smoothness of the solution, the smoothness of the data, and the numerical errors/order-of-accuracy of our scheme. We perform a detailed numerical study to quantify these errors experimentally. We start with a comparison of the numerical schemes for solving the hyperbolic PDEs that appear in the optimality system and the Newton step (see \secref{s:experiments:pdesolver}). The second set of experiments analyzes the effectiveness of our schemes for preconditioning the reduced space KKT system (see \secref{s:experiments:precond}).

\begin{figure}
\centering
\includegraphics[width=0.6\textwidth]
{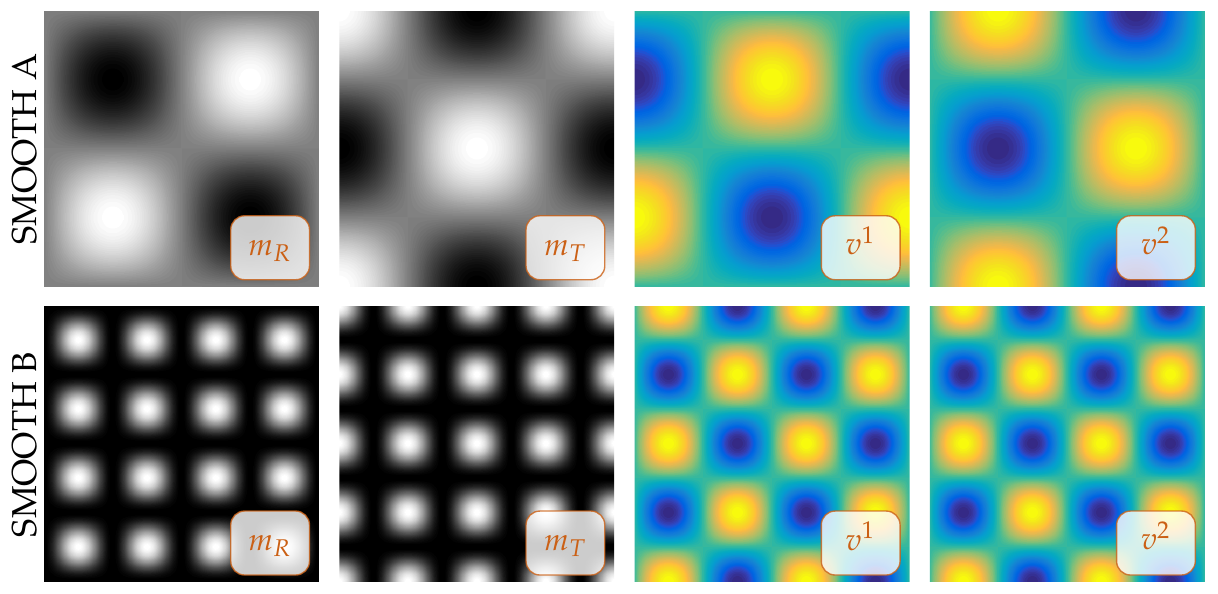}
\caption{Synthetic test problems. From left to right: reference image $m_R$; template image $m_T$; $v^1$ component of velocity field; and $v^2$ component of velocity field. The intensity values of the images are in $[0,1]$. The magnitude of the velocity field is in $[-0.5,0.5]$ (top row) and $[-1,1]$ (bottom row), respectively.
\label{f:synthetic-data}}
\end{figure}

\begin{figure}
\centering
\includegraphics[width=0.9\textwidth]
{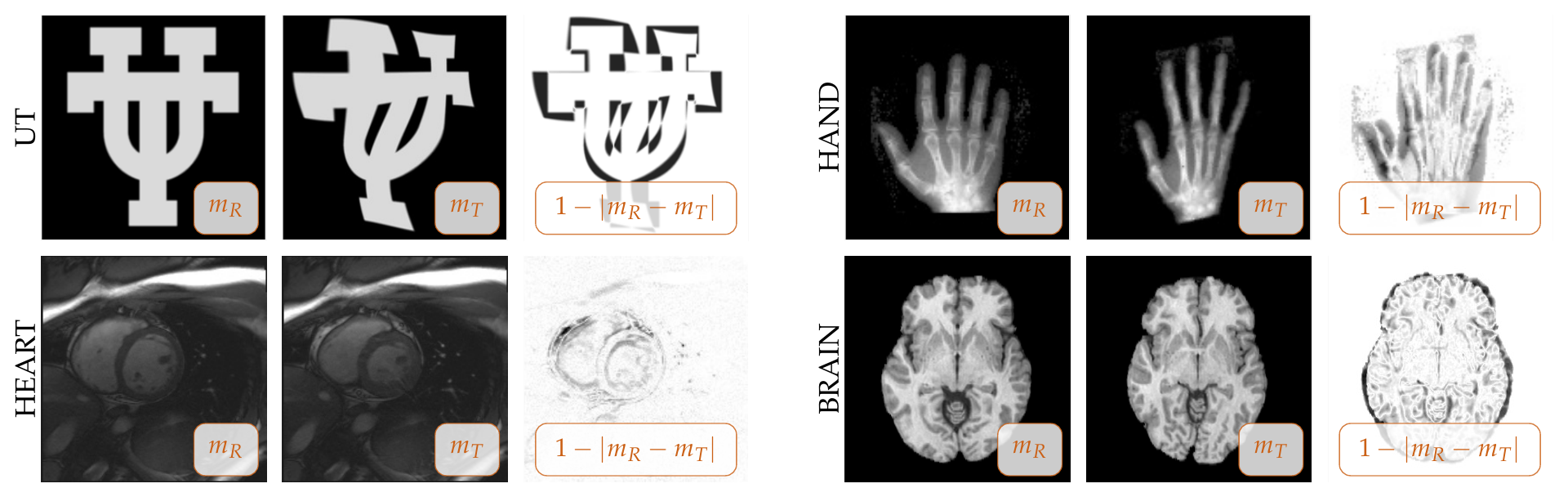}
\caption{Registration problems. Top left: UT images (synthetic problem); top right: HAND images~\cite{Amit:1994a,Modersitzki:2009a}; bottom left: HEART images; bottom right: BRAIN images~\cite{Christensen:2006a}. The intensity values for these images are normalized to $[0,1]$. We provide (from left to right for each set of images) the reference image $m_R$, the template image $m_T$, and the residual differences between these images prior to registration.
\label{f:registration-problems}}
\end{figure}

All experiments are carried out for $d=2$ using Matlab R2013a on a Linux cluster with Intel Xeon X5650 Westmere EP 6-core processors at 2.67GHz with 24GB DDR3-1333 memory. We illustrate the synthetic and real world data used for the experiments in~\figref{f:synthetic-data} and~\figref{f:registration-problems}, respectively.\footnote{The HAND images in \figref{f:registration-problems} are taken from~\cite{Modersitzki:2009a}. The BRAIN images in \figref{f:registration-problems} are taken from the `'Nonrigid Image Registration Evaluation Project`' (NIREP) available at \href{http://nirep.org}{\tt http://nirep.org} (data sets {\tt na01} and {\tt na02})~\cite{Christensen:2006a}.}

\subsection{Hyperbolic PDE solver}
\label{s:experiments:pdesolver}

We study the performance of the time integrators for the hyperbolic transport equations. We only consider the problems SMOOTH A and SMOOTH B in \figref{f:synthetic-data} as these are constructed to be initially resolved on the considered grids $\Omega^h$. This allows us to study grid convergence without mixing in any additional problems due to potential sharp transitions in the intensity values of the image data. We will see that these simple test cases can already break standard numerical schemes.

\subsubsection{Self-Convergence: State and Adjoint Equation}

\ipoint{Purpose} To study the numerical stability and accuracy of the considered schemes for integrating the hyperbolic transport equations that appear in our optimality system.

\ipoint{Setup} We study the self-convergence of the considered numerical time integrators. We consider the RK2 scheme (pseudospectral discretization in space), the stabilized RK2A scheme (pseudospectral discretization in space), and the SL method (cubic interpolation combined with a pseudospectral discretization; see \secref{s:timeintegration} for details). We test these schemes for the synthetic problems SMOOTH A and SMOOTH B in \figref{f:synthetic-data}. We consider the state and the adjoint equation. We compute the relative $\ell^2$-error between the solution of the transport equations (state equation~\eqref{e:state-pde-elim} and adjoint equation~\eqref{e:adj-pde-elim}) obtained on a spatial grid of size $\vect{n}_x$ and the solution obtained on a spatial grid of size $\tilde{\vect{n}}_x = 2\vect{n}_x$. We compute this error in the Fourier domain; formally, the error is given by
\[
\|\delta u^h\|_{\text{rel}}
\defeq
{\|\mat{M}[\mat{W}^{-1}u^h]_{\vect{n}_x} - [\mat{W}^{-1}u^h]_{\tilde{\vect{n}}_x}\|_2}
/{\|[\mat{W}^{-1}u^h]_{\tilde{\vect{n}}_x}\|_2}
\]

\noindent for a given numerical solution $u^h$. Here, $[\,\cdot\,]_{\vect{n}}$ indicates that the data is represented on a grid of size $\vect{n}$; $\mat{M}$ is a prolongation operator that maps the data from a grid of size $\vect{n}_x$ to a grid of size $\tilde{\vect{n}}_x$; and $\mat{W}^{-1}$ represents the forward Fourier operator. We use a CFL number of 0.2 to compute the number of time steps $n_t$ for the RK2 and the RK2A method. For the SL method we use the CFL numbers 0.2, 1, and 5. We expect the error to tend to zero for an increasing number of discretization points.

\begin{figure}
\centering
\includegraphics[width=0.6\textwidth]
{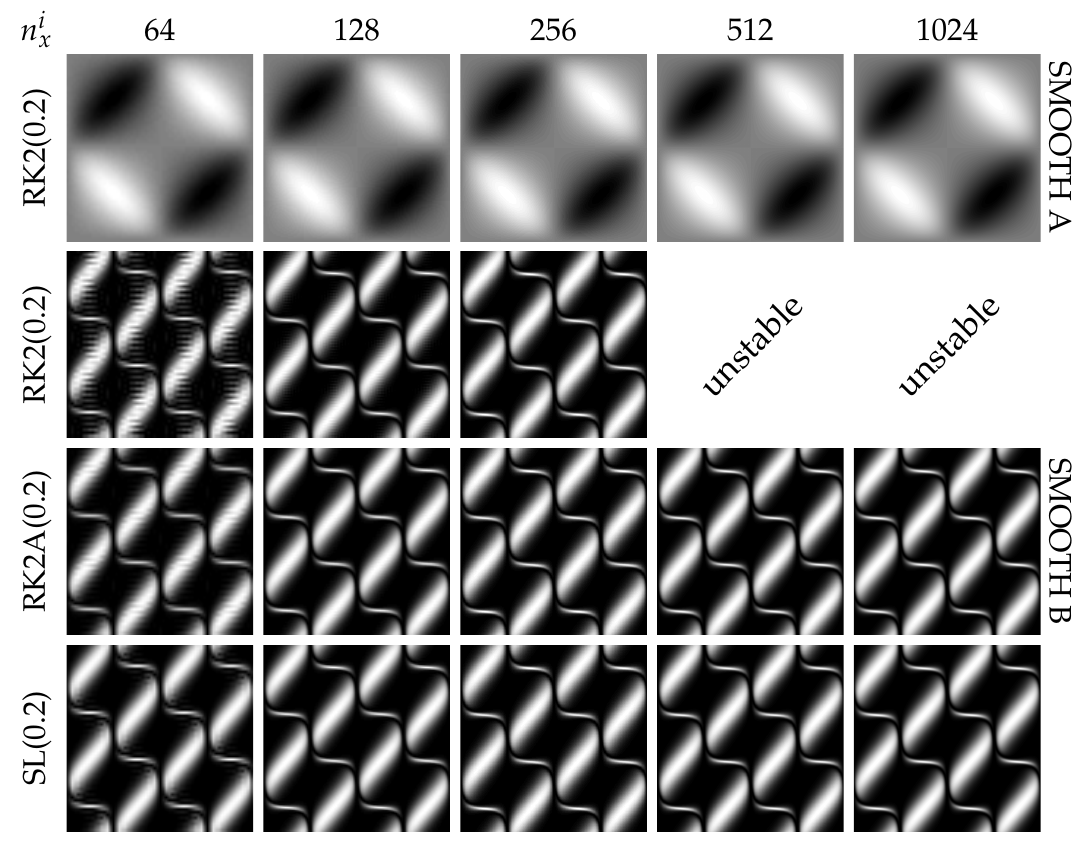}
\caption{Self-convergence for the forward solver. We illustrate solutions of the forward problem (state equation; see \eqref{e:state-pde-elim}) for the synthetic test problems in \figref{f:synthetic-data} (top row: SMOOTH A; bottom rows: SMOOTH B). We report results for different grid sizes $\vect{n}_x = (n^1_x,n^2_x)^\T$. We use the same number of time steps (CFL number of 0.2) for all PDE solvers.
\label{f:selfconvergence-fwdsolvers-imagecollection}}
\end{figure}

\ipoint{Results} We report results for the self-convergence of our numerical schemes in \tabref{t:selfconvergence-hyperbolic-pdesolvers}. We illustrate a subset of these results in \figref{f:selfconvergence-fwdsolvers-imagecollection}.

\ipoint{Observations} The most important observations are that \bipa\item our SL scheme delivers an accuracy that is at the order of the RK2A and the RK2 scheme with a speed up of one order of magnitude,\footnote{We have also tested an implementation of the SL method that delivers more accurate solutions (less numerical diffusion) at the expense of a significant increase in time to solution. In this scheme we upsampled the data to a grid of size $2\vect{n}_x$ whenever we had to interpolate. The associated gain in accuracy did not justify the increase in computational cost.} and \item that our standard RK2 scheme can become unstable if we combine it with a spectral discretization---even for smooth initial data and a smooth velocity field (\runref{31}, \runref{36}, \runref{71}, and \runref{76} in \tabref{t:selfconvergence-hyperbolic-pdesolvers})\eipa. This instability is a consequence of the absence of numerical diffusion; it is completely unrelated to the CFL condition. The RK2A and the SL method remain stable across all considered test cases with a similar performance. The rate of convergence for the RK2 and the RK2A scheme are excellent; we expect second order convergence in time and spectral convergence in space (this has been verified; results not reported here). The self-convergence for the SL(0.2) method is at the order, but overall slightly better, than the one observed for the RK2 and the RK2A scheme. The error for the self-convergence increases by one order of magnitude if we increase the CFL number for the SL method to 1 or 5, respectively. Switching from test problem SMOOTH A to SMOOTH B the self-convergence deteriorates for both methods. We can observe that we can not fully resolve the problem SMOOTH B if we solve the equations on a spatial grid with less than 128 nodes along each spatial direction; the errors range between $\bigO(\num{1E-1})$ and $\bigO(\num{1E-2})$ (\runref{21} through \runref{25} for the state equation and \runref{61} through \runref{70} for the adjoint equation; see also \figref{f:selfconvergence-fwdsolvers-imagecollection}). Notice that we can fully resolve the initial data and the velocity field for smaller grid sizes.

We can also observe that we loose about one order of magnitude in the rate of convergence if we switch from the state to the adjoint equation---even for the mild case SMOOTH A. This observation is consistent across all solvers. This demonstrates that the adjoint equation is in general more difficult to solve than the state equation. This can be attributed to the fact that the adjoint equation is a transport equation for the residual; the residual has, in general, less regularity than the original images (see also~\cite{Vialard:2012a}).

As for the time to solution we can observe that the SL(0.2) scheme delivers a performance that is at the order of the RK2A(0.2) scheme (slightly worse). We have to switch to use a CFL number of 1 to be competitive with the RK2(0.2) scheme. For a CFL number of 5 the SL scheme outperforms the RK2 and RK2A scheme by about one order of magnitude in terms of time to solution. Intuitively, one would expect that the SL method delivers much more pronounced speedup due to the unconditional stability. However, the discrepancy is due to the fact that we essentially replace a large number of highly optimized FFT operations with cubic spline interpolation operations. We report estimates for computational complexity in terms of FFTs and IPs in \tabref{t:computational-complexity-theoretical} in \secref{s:computational-complexity}. We can see in \tabref{t:timings-fft-vs-interpolation} in \secref{s:computational-complexity} that the differences in CPU time between these two operations are significant.

\ipoint{Conclusions} We can not guarantee convergence to a valid solution if we use a standard RK2 scheme in combination with a spectral discretization, even for smooth initial data; we have to use more sophisticated schemes. We provide two alternatives: a stabilized RK2 scheme (RK2A) and an SL scheme (see \secref{s:timeintegration} for details). Both schemes remained stable across all experiments. The SL scheme delivers a performance that is very similar to the RK2A scheme with a speedup of one order of magnitude---even for our non-optimized implementation.\footnote{Matlab's FFT library is based on the highly optimized FFTW library (see \url{http://www.fftw.org}; \cite{Frigo:2005a}). For interpolation we use Matlab's built in {\tt interp2} routine (Matlab R2013a; we report timings in \tabref{t:timings-fft-vs-interpolation} in \secref{s:computational-complexity}). We expect an additional speedup if we switch to an optimized, three-dimensional {\tt C++} implementation, something we will investigate in future work.}

\subsubsection{Convergence to RK2A}

\ipoint{Purpose} To assess \bipa\item the convergence of the SL method to the solution of the RK2A scheme and by that \item the numerical errors that might affect the overall convergence of our Newton--Krylov solver\eipa.

\ipoint{Setup} We assess the convergence of the SL method to a solution computed on the basis of the RK2A scheme for the state and the adjoint equation \eqref{e:state-pde-elim} and~\eqref{e:adj-pde-elim}, respectively. Based on our past experiments (see~\cite{Mang:2015a,Mang:2016a}) we assume that the solution of the RK2A scheme is a silver standard. We compute the reference solution on a grid of size $\vect{n}_x = (512,512)^\T$ with a CFL number of 0.2 (RK2A(0.2)). Likewise to the former experiment we compute the discrepancy between the numerical solutions in the Fourier domain; i.e., we report relative errors $\|\delta m^h_1\|_{\text{rel}}$ and $\|\delta \lambda^h_0\|_{\text{rel}}$, where we have $\tilde{\vect{n}}_x = (512,512)^\T$ for the RK2A(0.2) reference solution. We report results for different discretization levels (varying number of grid points $\vect{n}_x$ and $n_t$; the CFL numbers for the SL method are 0.2, 1, 2, 5, and 10). As a reference, we also compute errors for the RK2A scheme for a CFL number of 0.2. We also compute convergence errors for the gradient; the setup is the same as for the experiment for the adjoint and state equation.

\ipoint{Results} We report the relative error between the solution computed based on our SL formulation and the RK2A scheme in \tabref{t:convergence-SL-to-RK2A-transporteqs}. The error estimates for the reduced gradient can be found in \tabref{t:convergence-SL-to-RK2A-gradient}.

\begin{table}
\caption{Convergence of the reduced gradient computed via the SL method to the gradient computed via the RK2A method. We evaluate the reference gradient on a grid of size $\vect{n}_x=(512,512)^\T$ via the RK2A method with a CFL number of 0.2. For the SL method the reduced gradient is computed on a grid of size $\vect{n}_x=(256,256)^\T$ and $\vect{n}_x=(512,512)^\T$ with a varying number of time steps $n_t$. We report the CFL number $c$, the associated number of time steps $n_t$, the relative $\ell^2$-error between numerical approximations to the reduced gradient $\vect{g}^h$, and the wall-clock time for the evaluation of $\vect{g}^h$. We consider the test problems SMOOTH A and SMOOTH B in \figref{f:synthetic-data} as input data. As a reference, we also provide relative errors for the RK2A scheme.
\label{t:convergence-SL-to-RK2A-gradient}}
\centering
\resetrunid
\begin{scriptsize}
\begin{tabular}{rrrrRrrrrrRrrr}
\toprule
  \mcol{8}{SMOOTH A} & \mcol{6}{SMOOTH B}
\\\midrule
$n_x^i$ & $c$ & run & $n_t$ & SL & time & RK2A & time & run & $n_t$ & SL & time & RK2A & time
\\\midrule
  256
& 10
& \runid &   3 & \num{2.538743e-03} & \num{4.511290e-01} & ---                & ---
& \runid &   9 & \num{2.282562e-02} & \num{1.292376e+00} & ---                & ---
\\
& 5
& \runid &   5 & \num{9.079548e-04} & \num{1.019461e+00} & ---                & ---
& \runid &  17 & \num{2.216633e-02} & \num{2.186887e+00} & ---                & ---
\\
& 2
& \runid &  11 & \num{2.191465e-04} & \num{1.767848e+00} & ---                & ---
& \runid &  41 & \num{2.206116e-02} & \num{3.869897e+00} & ---                & ---
\\
& 1
& \runid &  21 & \num{1.306122e-04} & \num{2.589129e+00} & ---                & ---
& \runid &  82 & \num{2.204569e-02} & \num{9.554244e+00} & ---                & ---
\\
& 0.2
& \runid & 102 & \num{1.209628e-04} & \num{1.170669e+01} & \num{1.209824e-04} & \num{8.565822e+00}
& \runid & 408 & \num{2.203461e-02} & \num{2.998994e+01} & \num{2.192176e-02} & \num{2.741565e+01}
\\\midrule
  512
& 10
& \runid &   5 & \num{8.998493e-04} & \num{3.662705e+00} & ---                & ---
& \runid &  17 & \num{1.415679e-03} & \num{7.587263e+00} & ---                & ---
\\
& 5
& \runid &   9 & \num{2.741060e-04} & \num{3.738178e+00} & ---                & ---
& \runid &  33 & \num{3.936790e-04} & \num{1.189113e+01} & ---                & ---
\\
& 2
& \runid &  21 & \num{4.927827e-05} & \num{8.719499e+00} & ---                & ---
& \runid &  82 & \num{7.890460e-05} & \num{2.855425e+01} & ---                & ---
\\
& 1
& \runid &  41 & \num{1.233117e-05} & \num{1.493975e+01} & ---                & ---
& \runid & 163 & \num{3.599177e-05} & \num{6.133237e+01} & ---                & ---
\\
& 0.2
& \runid & 204 & \num{4.871694e-07} & \num{7.319492e+01} & 0                  & \num{4.962631e+01}
& \runid & 815 & \num{2.902154e-05} & \num{2.654561e+02} & 0                  & \num{2.036633e+02}
\\\bottomrule
\end{tabular}
\end{scriptsize}
\end{table}

\ipoint{Observations} The most important observation is that the SL scheme converges to the RK2A(0.2) reference solution with a similar rate than the RK2A scheme itself.  The SL scheme delivers an equivalent or even better rate of convergence than the RK2A scheme for a CFL number of 0.2. We lose one to two digits if we switch to higher CFL numbers; this loss in accuracy might still be acceptable for our Newton--Krylov solver to converge to almost identical solutions, something we will investigate below. Likewise to the former experiment we can again observe that the error for the adjoint equation are overall about one order of magnitude larger than those obtained for the state equation; this observation is again consistent for both schemes---the RK2A scheme and the SL scheme (see for instance \runref{41} and \runref{45}; and \runref{42} and \runref{46} in \tabref{t:convergence-SL-to-RK2A-transporteqs}).

\ipoint{Conclusions} Our SL scheme behaves very similar than the RK2A scheme with the benefit of an orders of magnitude reduction in computational work load due to the unconditional stability. We expect significant savings, especially for evaluating the Hessian, as accuracy requirements for the Hessian and its preconditioner are less significant than those for the reduced gradient for our Newton--Krylov solver to still converge. If high accuracy solutions are required, we can simply increase the number of time points to match the accuracy obtained for the RK2A scheme at the expense of an increase in CPU time.

\subsubsection{Adjoint Error}
\label{s:experiments:adj-err}

\ipoint{Purpose} To assess the numerical errors of the discretized forward and adjoint operator.

\ipoint{Setup} We solve the state equation~\eqref{e:state-pde-elim} and the adjoint equation~\eqref{e:adj-pde-elim} on a grid of size $\vect{n}_x = (256,256)^\T$ for a varying number of time points $n_t$. We consider the problem SMOOTH A in \figref{f:synthetic-data} to setup the equations. We report the relative error between the discretized forward operator $\D{C}^h$ and the discretized adjoint operator $(\D{C}^h)^\T$: $\delta_{\text{ADJ}}\defeq{|\langle\D{C}^h m_0^h, \D{C}^h m_0^h\rangle-\langle(\D{C}^h)^{\T}\D{C}^h m_0^h, m_0^h\rangle|}/{|\langle\D{C}^h m_0^h, \D{C}^h m_0^h\rangle|}$. The continuous forward operator is self-adjoint, i.e., $\D{C} = \D{C}^\T$; the error should tend to zero if our numerical scheme preserves this property.\footnote{Our solver is based on an optimize-then-discretize approach (see~\secref{s:solver}). We can not guarantee that the properties of the continuous operators of our constrained formulation and its variations are preserved after discretization. In a discretize-then-optimize approach the discretization is differentiated, which will result in consistent operators; we refer to~\cite[pages~57ff.]{Gunzburger:2003a} for a more detailed discussion on the pros and cons; we also discuss this in \secref{s:solver}.}

\ipoint{Results} We report the relative adjoint errors in \tabref{t:adjoint-error}.

\begin{table}
\caption{Relative adjoint error $\delta_{\text{ADJ}}$ (see text for details) for a grid size of $\vect{n}_x = (256,256)^\T$ and a varying number of time steps $n_t$. We consider the test problem SMOOTH A in \figref{f:synthetic-data}. We report the CFL number $c$, the associated number of time steps $n_t$, and the relative errors for the SL and the RK2A scheme.
\label{t:adjoint-error}}
\centering
\begin{scriptsize}
\begin{tabular}{rrRr}
\toprule
$c$ & $n_t$ & SL                 & RK2A               \\\midrule
 10 &     3 & \num{3.280379e-03} & ---                \\
  5 &     5 & \num{1.262295e-03} & ---                \\
  2 &    11 & \num{2.754930e-04} & ---                \\
  1 &    21 & \num{7.723612e-05} & ---                \\
0.2 &   102 & \num{3.296296e-06} & \num{1.244068e-16} \\\bottomrule
\end{tabular}
\end{scriptsize}
\end{table}

\ipoint{Observations} The most important observation is that the RK2A scheme is self-adjoint (up to machine precision) whereas the error for the SL method ranges between $\bigO(\num{1E-6})$ and $\bigO(\num{1E-3})$ as a function of $n_t$. If we solve the problem with a CFL number of 2 or smaller, the adjoint error is below or at the order of the accuracy we typically solve the inverse problem with in practical applications (relative change of the gradient of $\num{1E-2}$ or $\num{1E-3}$ and an absolute tolerance for the $\ell^\infty$-norm of the reduced gradient of $\num{1E-5}$).

\ipoint{Conclusions} Our SL scheme is not self-adjoint. The numerical errors are acceptable for the tolerances we use in practical applications---even for moderate CFL numbers. If we intend to solve the problem with a higher accuracy, we might have to either use a larger number of time steps or switch to the RK2A scheme to guarantee convergence. We already note that we have not observed any problems in terms of the convergence (failure to converge) nor the necessity for any additional line search steps in our solver, even if we considered  a CFL number of 10.

\subsection{Preconditioner}
\label{s:experiments:precond}

Next, we analyze the performance of our preconditioners (see \secref{s:newtonkrylovmethod}).

\subsubsection{Eigenvalue Estimation}
\label{s:experiments:eigenvalest}

We need to estimate the extremal eigenvalues of $\mat{P}_{\text{2L}}$ if we use the CHEB method to compute the action of its inverse. This estimation results in a significant amount of computational work if we have to do it frequently (about 30 matvecs for the estimation of $e_{\max}$). We estimate the smallest eigenvalue based on an analytical approximation; we estimate the largest eigenvalue numerically (see~\secref{s:precond} for details).

\ipoint{Purpose} To assess if the estimates for the largest eigenvalue vary significantly during the course of on inverse solve.

\ipoint{Setup} We solve the inverse problem for different sets of images; we consider the UT, HAND, HEART, and BRAIN images in~\figref{f:registration-problems}. We terminate the inversion if the gradient is reduced by three orders of magnitude or if $\|\vect{g}^h_k\|_{\infty} \leq \num{1E-5}$, $k=0,1,2,\ldots$. We estimate the largest eigenvalue every time the preconditioner is applied. We consider a compressible diffeomorphism ($H^2$-regularization). The solution is computed using a GN approximation. We report results for different regularization weights $\beta_v$.

\ipoint{Results} We summarize the estimates for the largest eigenvalue $e_{\max}$ in \tabref{t:eigenvalue-estimates-invsolve}.

\begin{table}
\caption{Estimates for the largest eigenvalue $e_{\max}$ of $\mat{P}_{\text{2L}}$ during the course of the inversion. We limit this experiment to a compressible diffeomorphism ($H^2$-regularization). We report results for the UT ($256\times256$), the HAND ($128\times128$), the HEART ($192\times192$), and the BRAIN ($256\times300$) images (see \figref{f:registration-problems}). We terminate the inversion if the relative change of the gradient is equal or smaller than three orders of magnitude or if the $\ell^\infty$-norm of the reduced gradient is smaller or equal to $\num{1E-5}$. We consider different regularization weights $\beta_v$. We estimate $e_{\max}$ every time we apply the preconditioner $\mat{P}_{\text{2L}}$. We report the initial estimate for $e_{\max}$ (zero velocity field), and the $\min$, $\operatorname{mean}$, and $\max$ values of the estimates computed during the course of the entire inversion.
\label{t:eigenvalue-estimates-invsolve}}
\centering
\begin{scriptsize}
\resetrunid
\begin{tabular}{lrlLlll}
\toprule
       & run    & $\beta_v$    & $e_{\max,0}$       & $\min$             & $\max$             & $\operatorname{mean}$ \\\midrule
UT     & \runid & $\num{1E-1}$ & \num{1.065406e+02} & \num{9.260466e+01} & \num{1.065406e+02} & \num{9.342118e+01} \\
       & \runid & $\num{1E-2}$ & \num{1.056406e+03} & \num{8.186378e+02} & \num{1.056406e+03} & \num{8.594239e+02} \\
       & \runid & $\num{1E-3}$ & \num{1.055506e+04} & \num{7.879599e+03} & \num{1.055506e+04} & \num{8.362990e+03} \\\midrule
HAND   & \runid & $\num{1E-1}$ & \num{2.729451e+01} & \num{2.501282e+01} & \num{2.729451e+01} & \num{2.535642e+01} \\
       & \runid & $\num{1E-2}$ & \num{2.639451e+02} & \num{2.203823e+02} & \num{2.639451e+02} & \num{2.239624e+02} \\
       & \runid & $\num{1E-3}$ & \num{2.630451e+03} & \num{2.133649e+03} & \num{2.630451e+03} & \num{2.161678e+03} \\\midrule
HEART  & \runid & $\num{1E-1}$ & \num{4.224790e+01} & \num{4.224790e+01} & \num{4.230641e+01} & \num{4.229045e+01} \\
       & \runid & $\num{1E-2}$ & \num{4.134790e+02} & \num{4.134790e+02} & \num{4.138565e+02} & \num{4.138241e+02} \\
       & \runid & $\num{1E-3}$ & \num{4.125790e+03} & \num{4.123857e+03} & \num{4.126295e+03} & \num{4.125989e+03} \\\midrule
BRAIN  & \runid & $\num{1E-1}$ & \num{2.970803e+01} & \num{2.918645e+01} & \num{2.970803e+01} & \num{2.937367e+01} \\
       & \runid & $\num{1E-2}$ & \num{2.880803e+02} & \num{2.816701e+02} & \num{2.904633e+02} & \num{2.893691e+02} \\
       & \runid & $\num{1E-3}$ & \num{2.871803e+03} & \num{2.807852e+03} & \num{2.871803e+03} & \num{2.825020e+03} \\\bottomrule
\end{tabular}
\end{scriptsize}
\end{table}

\ipoint{Observations} The most important observation is that the estimates for the largest eigenvalue do not vary significantly during the course of the iterations for most of the considered test cases. We have verified this for different reference and template images and as such for varying velocity fields. Our results suggest that we might have to only estimate the eigenvalues once for the initial guess---a zero velocity field. The costs for applying the Hessian for a zero velocity field are small---several expressions in~\eqref{e:hessian-matvec},~\eqref{e:inc-state-pde-elim}, and~\eqref{e:inc-adj-pde-elim} drop or are constant. Our results suggest that the changes in the eigenvalues are a function of the changes in the magnitude of the velocity field $\vect{v}$, i.e., the amount of expected deformation between the images. That is, we have only subtle residual differences and a small deformation in case of the HEART images; the estimated eigenvalues are almost constant. For the HAND and the UT images the deformations and the residual differences are larger; the changes in the estimates for the largest eigenvalue are more pronounced. Another important observation is that the most significant changes occur during the first few outer iterations. Once we are close to the solution of our problem, the eigenvalues are almost constant.

Overall, these results suggest that we can limit the estimation of the eigenvalues to the first iteration, or---if we observe a deterioration in the performance of our preconditioner---re-estimate the eigenvalues and for the subsequent solves again keep them fixed. We can observe that we have to estimate the eigenvalues only once for a given set of images; changes in the regularization parameter can simply be accounted for by rescaling these eigenvalue estimates. This is in accordance with our theoretical understanding of how changes in the regularization parameter affect the spectrum of the Hessian operator.

\ipoint{Conclusions} The estimates for the largest eigenvalue do not vary significantly during the course of the inversion for the considered test problems. We can estimate the eigenvalues efficiently during the first iteration (zero initial guess) and potentially use this estimate throughout the entire inversion.

\subsubsection{Convergence: KKT Solve}

\ipoint{Purpose} To assess the rate of convergence of the KKT solve for the different schemes to precondition the reduced space Hessian.

\ipoint{Setup} We consider three sets of images, the test problem SMOOTH A in \figref{f:synthetic-data}, and the BRAIN and the HAND images in \figref{f:registration-problems}. We solve the forward problem to setup a synthetic test problem based on the velocity field $\vect{v}^{h\star}$ of problem SMOOTH A, i.e., we transport $m_R$ to obtain a synthetic template image $m_T$. We consider a GN approximation to $\D{H}^h$. We study three schemes to precondition the KKT system: \bipa\item the regularization preconditioner $\mat{P}_{\text{REG}}$, \item the nested preconditioner $\mat{P}_{\text{2L}}$ the inverse action of which we compute using a PCG method, and \item the nested preconditioner $\mat{P}_{\text{2L}}$ the inverse action of which we compute based on a CHEB method\eipa. If we use a PCG method to invert the preconditioner, we have to use a higher accuracy than the one we use to solve the KKT system. We increase the accuracy by one order of magnitude; we refer to this solver as PCG($\num{1E-1}$). For the CHEB method we can use a fixed number of iterations; we have tested 5, 10, and 20 iterations. We observed an overall good performance for 10 iterations. We refer to this strategy as CHEB(10).

We perform two experiments: In the first experiment we use a true solution $\vect{\tilde{v}}^{h\star} = -0.5\vect{v}^{h\star}$ and apply the Hessian operator to generate a synthetic right hand side $\vect{\tilde{b}}\mspace{-2mu}{}^{\,h}$. We solve the KKT system $\D{H}^h\vect{\tilde{v}}^h = \vect{\tilde{b}}\mspace{-2mu}{}^{\,h}$ with a zero initial guess for $\vect{\tilde{v}}^h$ using a PCG method with a tolerance of \num{1E-12}. We compute the (relative) $\ell^2$-norm of the difference between $\vect{\tilde{v}}^h$ and $\vect{\tilde{v}}^{h\star}$ to assess if our schemes converge to the true solution with the same accuracy. We set up the KKT system based on the test problem SMOOTH A in \figref{f:synthetic-data}.

For the second experiment, we evaluate the reduced gradient $\vect{g}^h$ at the true solution $\vect{v}^{h\star}$ and solve the system $\D{H}^h\vect{\tilde{v}}^h = -\vect{g}^h$ with a zero initial guess for $\vect{\tilde{v}}^h$. We solve the system using a PCG method with a tolerance of $\num{1E-6}$. We consider a compressible diffeomorphism ($H^2$-regularization). We use the RK2A scheme with a CFL number of 0.2 for the regularization preconditioner and the SL scheme with a CFL number of five for the two-level preconditioner. We report results for the test problems SMOOTH A, BRAIN, and HAND. We consider different spatial resolution levels (grid convergence) and different choices for the regularization parameter $\beta_v$. An ideal preconditioner is mesh-independent and delivers the same rate of convergence irrespective of the choice of the regularization weight.

\begin{figure}
\centering
\includegraphics[width=0.9\textwidth]
{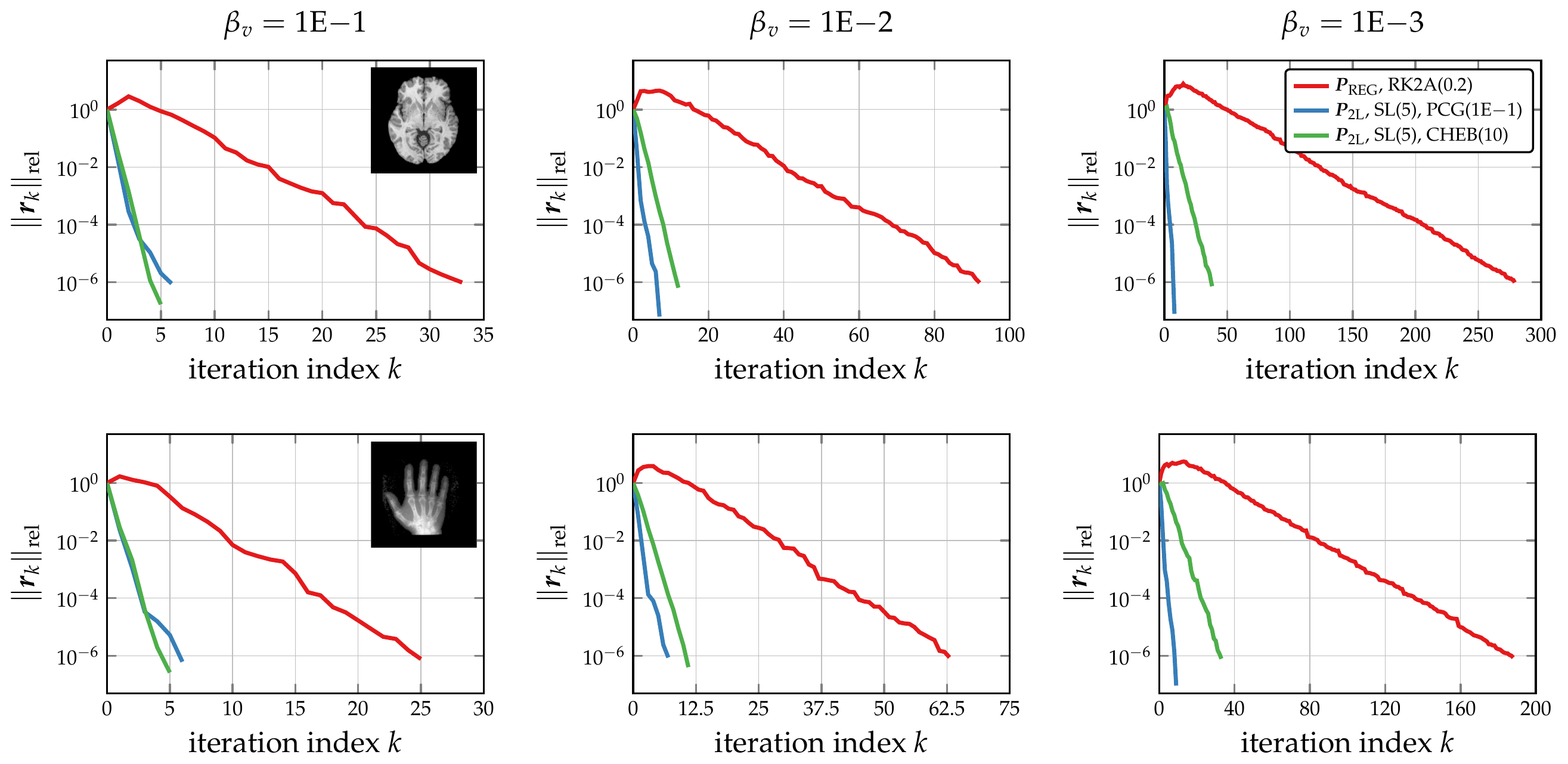}
\caption{Convergence results for different strategies to precondition the reduced space KKT system. We report exemplary trends of the relative residual $\|\vect{r}_k\|_{\text{rel}}\defeq\|\vect{r}_k\|_2/\|\vect{r}_0\|_2$ with respect to the iteration number $k$. We report results for different images (top row: BRAIN; bottom row: HAND; grid size $(256,256)^\T$) with respect to varying regularization weights $\beta_v$ (left column: $\beta_v=\num{1E-1}$; middle column: $\beta_v=\num{1E-2}$; right column: $\beta_v=\num{1E-3}$). We solve the system at the true solution available for the considered synthetic test problems. We use an $H^2$-regularization model (compressible diffeomorphism). We use a PCG method with a tolerance of $\num{1E-6}$ to solve this system. We report results for the regularization preconditioner $\mat{P}_{\text{REG}}$ (red curve) and the nested preconditioner $\mat{P}_{\text{2L}}$. We use two solvers to invert the preconditioner: PCG(\num{1E-1}) (blue curve) and CHEB(10) (green curve). The results correspond to those reported in \tabref{t:preconditioners-kktsolve-convergence}.
\label{f:preconditioners-kktsolve-convergence-plots}}
\end{figure}

\ipoint{Results} We summarize the results for the first experiment in \tabref{t:preconditioners-kktsolve-solutionerror} and the results for the second part in \tabref{t:preconditioners-kktsolve-convergence}. We illustrate the convergence of a subset of the results reported in \tabref{t:preconditioners-kktsolve-convergence} in \figref{f:preconditioners-kktsolve-convergence-plots}.

\begin{table}
\caption{Error between the true solution $\vect{\tilde{v}}^{h\star}$ and the numerical solution $\vect{\tilde{v}}^{h}$ of the KKT system for different schemes to precondition the reduced space Hessian. We report the absolute and the relative $\ell^2$-error between $\vect{\tilde{v}}^{h\star}$ and $\vect{\tilde{v}}^{h}$. We report results for different preconditioners ($\mat{P}_{\text{REG}}$, $\mat{P}_{\text{2L}}$), different choices for the PDE solver (SL($c$) for different CFL numbers $c$ and RK2A(0.2)), and different choices for the method to solve for the action of the inverse of the preconditioner (CHEB(10) and PCG(\num{1E-1})). We solve for $\vect{\tilde{v}}^h$ using a PCG method with a tolerance of $\num{1E-12}$. We consider the test problem SMOOTH A as in \figref{f:synthetic-data} to set up the problem. We solve the system on a grid of size $256\times256$.
\label{t:preconditioners-kktsolve-solutionerror}}
\centering
\begin{scriptsize}
\resetrunid
\begin{tabular}{rlrrrR}
\toprule
  run  & PC & PDE solver         & PC solver & $\|\delta\|_2$  & $\|\delta\|_{2,\text{rel}}$ \\\midrule
\runid & $\mat{P}_{\text{REG}}$  & RK2A(0.2) & ---             & \num{7.492265e-13} & \num{4.138931e-14} \\
\runid & $\mat{P}_{\text{2L}}$   & SL(0.2)   & CHEB(10)        & \num{8.167023e-14} & \num{4.511685e-15} \\
\runid &                         & SL(1)     & PCG(\num{1E-1}) & \num{1.633962e-12} & \num{9.026450e-14} \\
\runid &                         & SL(1)     & CHEB(10)        & \num{1.225236e-13} & \num{6.768538e-15} \\
\runid &                         & SL(2)     & PCG(\num{1E-1}) & \num{1.139704e-12} & \num{6.296035e-14} \\
\runid &                         & SL(2)     & CHEB(10)        & \num{1.892358e-13} & \num{1.045390e-14} \\
\runid &                         & SL(5)     & PCG(\num{1E-1}) & \num{4.487734e-12} & \num{2.479146e-13} \\
\runid &                         & SL(5)     & CHEB(10)        & \num{5.218249e-13} & \num{2.882703e-14} \\
\runid &                         & SL(10)    & PCG(\num{1E-1}) & \num{2.055538e-12} & \num{1.135535e-13} \\
\runid &                         & SL(10)    & CHEB(10)        & \num{2.055538e-12} & \num{1.135535e-13}
\\\bottomrule
\end{tabular}
\end{scriptsize}
\end{table}

\begin{landscape}
\begin{table}
\caption{Convergence results for different strategies to precondition the reduced space KKT system. We report results for two preconditioners---our original preconditioner based on the regularization operator ($\mat{P}_{\text{REG}}$) and the proposed, nested preconditioner ($\mat{P}_{\text{2L}}$). We solve the reduced space KKT system via a PCG method with a tolerance of $\num{1E-6}$. We use two different solvers for the latter to invert the preconditioner---a PCG method with a tolerance that is $\num{1E-1}$ times the tolerance of the PCG method used to solve the reduced space KKT system (i.e., a tolerance of $\num{1E-7}$) and a CHEB method with a fixed number of 10 iterations. We consider a compressible diffeomorphism with an $H^2$-regularization model. We report results for different images (SMOOTH A, BRAIN, and HAND), for different regularization weights $\beta_v$, and a varying grid sizes $\vect{n}_x$ (grid convergence; number of unknowns $n=2n_x^1n_x^2$). We solve the reduced space KKT system at the true solution $\vect{v}^{h\star}$; the velocity field $\vect{v}^{h\star}$ corresponds to the test problem SMOOTH A. We consider the RK2A method with a CFL number of 0.2 for the regularization preconditioner and the SL scheme for the nested preconditioner with a CFL number of 5. We report ($i$) the number of PCG iterations until convergence, ($ii$) the time spent on the Hessian matvecs (in seconds), ($iii$) the percentage of that time spent on inverting the preconditioner (if applicable), and ($iv$) the speedup compared to our original preconditioner (regularization preconditioner in combination with the RK2A scheme).
\label{t:preconditioners-kktsolve-convergence}}
\centering
\resetrunid
\begin{tiny}
\begin{tabular}{rllrrrRllLrRllLrRllL}
\toprule
\mcol{10}{SMOOTH A} & \mcol{5}{HAND} & \mcol{5}{BRAIN}
\\\midrule
$n$ & $\beta_v$ & $\mat{P}$ & PDE solver & PC solver & run & iter & time & \% PC & speedup & run & iter & time & \% PC & speedup & run & iter & time & \% PC & speedup \\\midrule
  \num{8192}
& $\num{1E-1}$
& $\boldsymbol{P}_{\text{REG}}$
& RK2A(0.2)
& ---
& \runid
& 4 
& \num{2.668259e+00}
& ---
& ---
& \runid
& 19 
& \num{7.119794e+00}
& ---
& ---
& \runid
& 21 
& \num{8.389931e+00}
& ---
& ---
\\
&
& $\boldsymbol{P}_{\text{2L}}$
& SL(5)
& PCG(\num{1E-1})
& \runid
& 2 
& \num{3.082383e+00}
& \num[scientific-notation=fixed,retain-zero-exponent=false]{6.404230e+01}\%
& \num{8.656481e-01}
& \runid
& 7 
& \num{1.669183e+01}
& \num[scientific-notation=fixed,retain-zero-exponent=false]{9.042114e+01}\%
& \num{4.265436e-01}
& \runid
& 7 
& \num{1.660855e+01}
& \num[scientific-notation=fixed,retain-zero-exponent=false]{9.114893e+01}\%
& \num{5.051573e-01}
\\
&
& $\boldsymbol{P}_{\text{2L}}$
& SL(5)
& CHEB(10)
& \runid
& 2 
& \num{2.360904e+00}
& \num[scientific-notation=fixed,retain-zero-exponent=false]{5.581417e+01}\%
& \num{1.130185e+00}
& \runid
& 6 
& \num{6.190187e+00}
& \num[scientific-notation=fixed,retain-zero-exponent=false]{7.555583e+01}\%
& \num{1.150174e+00}
& \runid
& 7 
& \num{6.441593e+00}
& \num[scientific-notation=fixed,retain-zero-exponent=false]{7.569437e+01}\%
& \num{1.302462e+00}
\\
& $\num{1E-2}$
& $\boldsymbol{P}_{\text{REG}}$
& RK2A(0.2)
& ---
& \runid
& 6 
& \num{3.140557e+00}
& ---
& ---
& \runid
& 47 
& \num{1.594852e+01}
& ---
& ---
& \runid
& 53 
& \num{1.755570e+01}
& ---
& ---
\\
&
& $\boldsymbol{P}_{\text{2L}}$
& SL(5)
& PCG(\num{1E-1})
& \runid
& 2 
& \num{4.309290e+00}
& \num[scientific-notation=fixed,retain-zero-exponent=false]{7.355585e+01}\%
& \num{7.287876e-01}
& \runid
& 8 
& \num{4.601402e+01}
& \num[scientific-notation=fixed,retain-zero-exponent=false]{9.595138e+01}\%
& \num{3.466013e-01}
& \runid
& 8 
& \num{4.394630e+01}
& \num[scientific-notation=fixed,retain-zero-exponent=false]{9.610313e+01}\%
& \num{3.994807e-01}
\\
&
& $\boldsymbol{P}_{\text{2L}}$
& SL(5)
& CHEB(10)
& \runid
& 3 
& \num{3.193372e+00}
& \num[scientific-notation=fixed,retain-zero-exponent=false]{6.215236e+01}\%
& \num{9.834611e-01}
& \runid
& 8 
& \num{8.261571e+00}
& \num[scientific-notation=fixed,retain-zero-exponent=false]{7.929518e+01}\%
& \num{1.930446e+00}
& \runid
& 7 
& \num{6.979650e+00}
& \num[scientific-notation=fixed,retain-zero-exponent=false]{7.614633e+01}\%
& \num{2.515269e+00}
\\
& $\num{1E-3}$
& $\boldsymbol{P}_{\text{REG}}$
& RK2A(0.2)
& ---
& \runid
& 16 
& \num{6.106627e+00}
& ---
& ---
& \runid
& 138 
& \num{4.262210e+01}
& ---
& ---
& \runid
& 161 
& \num{4.956111e+01}
& ---
& ---
\\
&
& $\boldsymbol{P}_{\text{2L}}$
& SL(5)
& PCG(\num{1E-1})
& \runid
& 2 
& \num{7.772829e+00}
& \num[scientific-notation=fixed,retain-zero-exponent=false]{8.454532e+01}\%
& \num{7.856376e-01}
& \runid
& 10 
& \num{1.662692e+02}
& \num[scientific-notation=fixed,retain-zero-exponent=false]{9.869301e+01}\%
& \num{2.563439e-01}
& \runid
& 11 
& \num{2.235435e+02}
& \num[scientific-notation=fixed,retain-zero-exponent=false]{9.897738e+01}\%
& \num{2.217068e-01}
\\
&
& $\boldsymbol{P}_{\text{2L}}$
& SL(5)
& CHEB(10)
& \runid
& 6 
& \num{5.925945e+00}
& \num[scientific-notation=fixed,retain-zero-exponent=false]{7.531190e+01}\%
& \num{1.030490e+00}
& \runid
& 24 
& \num{2.160538e+01}
& \num[scientific-notation=fixed,retain-zero-exponent=false]{8.392701e+01}\%
& \num{1.972754e+00}
& \runid
& 22 
& \num{2.040142e+01}
& \num[scientific-notation=fixed,retain-zero-exponent=false]{8.293568e+01}\%
& \num{2.429297e+00}
\\\midrule
  \num{32768}
& $\num{1E-1}$
& $\boldsymbol{P}_{\text{REG}}$
& RK2A(0.2)
& ---
& \runid
& 4 
& \num{1.894343e+01}
& ---
& ---
& \runid
& 22 
& \num{4.308238e+01}
& ---
& ---
& \runid
& 26 
& \num{5.931509e+01}
& ---
& ---
\\
&
& $\boldsymbol{P}_{\text{2L}}$
& SL(5)
& PCG(\num{1E-1})
& \runid
& 2 
& \num{3.525265e+00}
& \num[scientific-notation=fixed,retain-zero-exponent=false]{4.719021e+01}\%
& \num{5.373619e+00}
& \runid
& 7 
& \num{2.734446e+01}
& \num[scientific-notation=fixed,retain-zero-exponent=false]{8.564008e+01}\%
& \num{1.575543e+00}
& \runid
& 7 
& \num{2.947105e+01}
& \num[scientific-notation=fixed,retain-zero-exponent=false]{8.816219e+01}\%
& \num{2.012656e+00}
\\
&
& $\boldsymbol{P}_{\text{2L}}$
& SL(5)
& CHEB(10)
& \runid
& 2 
& \num{4.155957e+00}
& \num[scientific-notation=fixed,retain-zero-exponent=false]{5.082264e+01}\%
& \num{4.558139e+00}
& \runid
& 6 
& \num{9.703700e+00}
& \num[scientific-notation=fixed,retain-zero-exponent=false]{6.491571e+01}\%
& \num{4.439789e+00}
& \runid
& 6 
& \num{9.042058e+00}
& \num[scientific-notation=fixed,retain-zero-exponent=false]{7.128061e+01}\%
& \num{6.559910e+00}
\\
& $\num{1E-2}$
& $\boldsymbol{P}_{\text{REG}}$
& RK2A(0.2)
& ---
& \runid
& 6 
& \num{2.350193e+01}
& ---
& ---
& \runid
& 54 
& \num{1.074328e+02}
& ---
& ---
& \runid
& 74 
& \num{1.564645e+02}
& ---
& ---
\\
&
& $\boldsymbol{P}_{\text{2L}}$
& SL(5)
& PCG(\num{1E-1})
& \runid
& 3 
& \num{7.136396e+00}
& \num[scientific-notation=fixed,retain-zero-exponent=false]{6.443852e+01}\%
& \num{3.293249e+00}
& \runid
& 7 
& \num{6.499782e+01}
& \num[scientific-notation=fixed,retain-zero-exponent=false]{9.488491e+01}\%
& \num{1.652868e+00}
& \runid
& 7 
& \num{7.240026e+01}
& \num[scientific-notation=fixed,retain-zero-exponent=false]{9.603799e+01}\%
& \num{2.161104e+00}
\\
&
& $\boldsymbol{P}_{\text{2L}}$
& SL(5)
& CHEB(10)
& \runid
& 3 
& \num{5.164036e+00}
& \num[scientific-notation=fixed,retain-zero-exponent=false]{5.589233e+01}\%
& \num{4.551078e+00}
& \runid
& 9 
& \num{1.262183e+01}
& \num[scientific-notation=fixed,retain-zero-exponent=false]{6.926115e+01}\%
& \num{8.511666e+00}
& \runid
& 10 
& \num{1.439546e+01}
& \num[scientific-notation=fixed,retain-zero-exponent=false]{6.992296e+01}\%
& \num{1.086902e+01}
\\
& $\num{1E-3}$
& $\boldsymbol{P}_{\text{REG}}$
& RK2A(0.2)
& ---
& \runid
& 16 
& \num{3.847602e+01}
& ---
& ---
& \runid
& 160 
& \num{3.452485e+02}
& ---
& ---
& \runid
& 224 
& \num{5.472086e+02}
& ---
& ---
\\
&
& $\boldsymbol{P}_{\text{2L}}$
& SL(5)
& PCG(\num{1E-1})
& \runid
& 3 
& \num{1.266749e+01}
& \num[scientific-notation=fixed,retain-zero-exponent=false]{7.787719e+01}\%
& \num{3.037383e+00}
& \runid
& 9 
& \num{2.531478e+02}
& \num[scientific-notation=fixed,retain-zero-exponent=false]{9.827132e+01}\%
& \num{1.363822e+00}
& \runid
& 10 
& \num{3.489000e+02}
& \num[scientific-notation=fixed,retain-zero-exponent=false]{9.882337e+01}\%
& \num{1.568382e+00}
\\
&
& $\boldsymbol{P}_{\text{2L}}$
& SL(5)
& CHEB(10)
& \runid
& 6 
& \num{8.505777e+00}
& \num[scientific-notation=fixed,retain-zero-exponent=false]{6.755726e+01}\%
& \num{4.523516e+00}
& \runid
& 27 
& \num{3.924277e+01}
& \num[scientific-notation=fixed,retain-zero-exponent=false]{7.390945e+01}\%
& \num{8.797761e+00}
& \runid
& 31 
& \num{4.098524e+01}
& \num[scientific-notation=fixed,retain-zero-exponent=false]{7.673541e+01}\%
& \num{1.335136e+01}
\\\midrule
  \num{131072}
& $\num{1E-1}$
& $\boldsymbol{P}_{\text{REG}}$
& RK2A(0.2)
& ---
& \runid
& 4 
& \num{6.620590e+01}
& ---
& ---
& \runid
& 25 
& \num{2.705618e+02}
& ---
& ---
& \runid
& 33 
& \num{3.519205e+02}
& ---
& ---
\\
&
& $\boldsymbol{P}_{\text{2L}}$
& SL(5)
& PCG(\num{1E-1})
& \runid
& 2 
& \num{1.242213e+01}
& \num[scientific-notation=fixed,retain-zero-exponent=false]{3.160107e+01}\%
& \num{5.329674e+00}
& \runid
& 6 
& \num{7.170931e+01}
& \num[scientific-notation=fixed,retain-zero-exponent=false]{7.858071e+01}\%
& \num{3.773036e+00}
& \runid
& 6 
& \num{1.016146e+02}
& \num[scientific-notation=fixed,retain-zero-exponent=false]{8.262478e+01}\%
& \num{3.463287e+00}
\\
&
& $\boldsymbol{P}_{\text{2L}}$
& SL(5)
& CHEB(10)
& \runid
& 2 
& \num{1.332895e+01}
& \num[scientific-notation=fixed,retain-zero-exponent=false]{4.648840e+01}\%
& \num{4.967075e+00}
& \runid
& 5 
& \num{2.584410e+01}
& \num[scientific-notation=fixed,retain-zero-exponent=false]{5.336051e+01}\%
& \num{1.046900e+01}
& \runid
& 5 
& \num{3.067467e+01}
& \num[scientific-notation=fixed,retain-zero-exponent=false]{5.695588e+01}\%
& \num{1.147267e+01}
\\
& $\num{1E-2}$
& $\boldsymbol{P}_{\text{REG}}$
& RK2A(0.2)
& ---
& \runid
& 6 
& \num{8.381351e+01}
& ---
& ---
& \runid
& 63 
& \num{6.606482e+02}
& ---
& ---
& \runid
& 92 
& \num{9.415809e+02}
& ---
& ---
\\
&
& $\boldsymbol{P}_{\text{2L}}$
& SL(5)
& PCG(\num{1E-1})
& \runid
& 2 
& \num{1.561411e+01}
& \num[scientific-notation=fixed,retain-zero-exponent=false]{3.961413e+01}\%
& \num{5.367806e+00}
& \runid
& 7 
& \num{2.165987e+02}
& \num[scientific-notation=fixed,retain-zero-exponent=false]{9.220416e+01}\%
& \num{3.050102e+00}
& \runid
& 7 
& \num{2.505603e+02}
& \num[scientific-notation=fixed,retain-zero-exponent=false]{9.407168e+01}\%
& \num{3.757901e+00}
\\
&
& $\boldsymbol{P}_{\text{2L}}$
& SL(5)
& CHEB(10)
& \runid
& 3 
& \num{2.064274e+01}
& \num[scientific-notation=fixed,retain-zero-exponent=false]{4.759811e+01}\%
& \num{4.060193e+00}
& \runid
& 11 
& \num{5.747926e+01}
& \num[scientific-notation=fixed,retain-zero-exponent=false]{5.738510e+01}\%
& \num{1.149368e+01}
& \runid
& 12 
& \num{5.859178e+01}
& \num[scientific-notation=fixed,retain-zero-exponent=false]{5.985408e+01}\%
& \num{1.607019e+01}
\\
& $\num{1E-3}$
& $\boldsymbol{P}_{\text{REG}}$
& RK2A(0.2)
& ---
& \runid
& 16 
& \num{1.660802e+02}
& ---
& ---
& \runid
& 188 
& \num{1.736728e+03}
& ---
& ---
& \runid
& 279 
& \num{2.293819e+03}
& ---
& ---
\\
&
& $\boldsymbol{P}_{\text{2L}}$
& SL(5)
& PCG(\num{1E-1})
& \runid
& 3 
& \num{3.082333e+01}
& \num[scientific-notation=fixed,retain-zero-exponent=false]{7.516487e+01}\%
& \num{5.388133e+00}
& \runid
& 9 
& \num{6.728488e+02}
& \num[scientific-notation=fixed,retain-zero-exponent=false]{9.763921e+01}\%
& \num{2.581156e+00}
& \runid
& 8 
& \num{9.542000e+02}
& \num[scientific-notation=fixed,retain-zero-exponent=false]{9.818625e+01}\%
& \num{2.403918e+00}
\\
&
& $\boldsymbol{P}_{\text{2L}}$
& SL(5)
& CHEB(10)
& \runid
& 6 
& \num{3.040049e+01}
& \num[scientific-notation=fixed,retain-zero-exponent=false]{5.777273e+01}\%
& \num{5.463076e+00}
& \runid
& 33 
& \num{1.366960e+02}
& \num[scientific-notation=fixed,retain-zero-exponent=false]{6.514839e+01}\%
& \num{1.270504e+01}
& \runid
& 38 
& \num{1.744101e+02}
& \num[scientific-notation=fixed,retain-zero-exponent=false]{6.302582e+01}\%
& \num{1.315187e+01}
\\\midrule
  \num{524288}
& $\num{1E-1}$
& $\boldsymbol{P}_{\text{REG}}$
& RK2A(0.2)
& ---
& \runid
& 4 
& \num{5.399723e+02}
& ---
& ---
& \runid
& 25 
& \num{1.974343e+03}
& ---
& ---
& \runid
& 37 
& \num{2.684123e+03}
& ---
& ---
\\
&
& $\boldsymbol{P}_{\text{2L}}$
& SL(5)
& PCG(\num{1E-1})
& \runid
& 2 
& \num{6.196909e+01}
& \num[scientific-notation=fixed,retain-zero-exponent=false]{2.900541e+01}\%
& \num{8.713575e+00}
& \runid
& 5 
& \num{3.539941e+02}
& \num[scientific-notation=fixed,retain-zero-exponent=false]{7.706754e+01}\%
& \num{5.577333e+00}
& \runid
& 5 
& \num{4.254987e+02}
& \num[scientific-notation=fixed,retain-zero-exponent=false]{8.331180e+01}\%
& \num{6.308181e+00}
\\
&
& $\boldsymbol{P}_{\text{2L}}$
& SL(5)
& CHEB(10)
& \runid
& 2 
& \num{7.150189e+01}
& \num[scientific-notation=fixed,retain-zero-exponent=false]{3.560023e+01}\%
& \num{7.551861e+00}
& \runid
& 4 
& \num{1.211438e+02}
& \num[scientific-notation=fixed,retain-zero-exponent=false]{4.649818e+01}\%
& \num{1.629752e+01}
& \runid
& 5 
& \num{1.600632e+02}
& \num[scientific-notation=fixed,retain-zero-exponent=false]{5.010754e+01}\%
& \num{1.676914e+01}
\\
& $\num{1E-2}$
& $\boldsymbol{P}_{\text{REG}}$
& RK2A(0.2)
& ---
& \runid
& 6 
& \num{6.487003e+02}
& ---
& ---
& \runid
& 67 
& \num{4.608890e+03}
& ---
& ---
& \runid
& 103 
& \num{7.124721e+03}
& ---
& ---
\\
&
& $\boldsymbol{P}_{\text{2L}}$
& SL(5)
& PCG(\num{1E-1})
& \runid
& 2 
& \num{7.122617e+01}
& \num[scientific-notation=fixed,retain-zero-exponent=false]{3.697354e+01}\%
& \num{9.107612e+00}
& \runid
& 6 
& \num{8.719053e+02}
& \num[scientific-notation=fixed,retain-zero-exponent=false]{9.068294e+01}\%
& \num{5.285998e+00}
& \runid
& 6 
& \num{1.203032e+03}
& \num[scientific-notation=fixed,retain-zero-exponent=false]{9.344625e+01}\%
& \num{5.922304e+00}
\\
&
& $\boldsymbol{P}_{\text{2L}}$
& SL(5)
& CHEB(10)
& \runid
& 3 
& \num{8.310675e+01}
& \num[scientific-notation=fixed,retain-zero-exponent=false]{3.817743e+01}\%
& \num{7.805627e+00}
& \runid
& 11 
& \num{2.961277e+02}
& \num[scientific-notation=fixed,retain-zero-exponent=false]{5.558788e+01}\%
& \num{1.556386e+01}
& \runid
& 14 
& \num{3.513835e+02}
& \num[scientific-notation=fixed,retain-zero-exponent=false]{5.599890e+01}\%
& \num{2.027620e+01}
\\
& $\num{1E-3}$
& $\boldsymbol{P}_{\text{REG}}$
& RK2A(0.2)
& ---
& \runid
& 16 
& \num{1.306679e+03}
& ---
& ---
& \runid
& 196 
& \num{1.307609e+04}
& ---
& ---
& \runid
& 310 
& \num{2.098427e+04}
& ---
& ---
\\
&
& $\boldsymbol{P}_{\text{2L}}$
& SL(5)
& PCG(\num{1E-1})
& \runid
& 2 
& \num{1.273722e+02}
& \num[scientific-notation=fixed,retain-zero-exponent=false]{5.983948e+01}\%
& \num{1.025875e+01}
& \runid
& 7 
& \num{2.926413e+03}
& \num[scientific-notation=fixed,retain-zero-exponent=false]{9.684914e+01}\%
& \num{4.468300e+00}
& \runid
& 7 
& \num{4.581961e+03}
& \num[scientific-notation=fixed,retain-zero-exponent=false]{9.799331e+01}\%
& \num{4.579757e+00}
\\
&
& $\boldsymbol{P}_{\text{2L}}$
& SL(5)
& CHEB(10)
& \runid
& 6 
& \num{1.681215e+02}
& \num[scientific-notation=fixed,retain-zero-exponent=false]{4.940214e+01}\%
& \num{7.772230e+00}
& \runid
& 35 
& \num{8.891574e+02}
& \num[scientific-notation=fixed,retain-zero-exponent=false]{6.045345e+01}\%
& \num{1.470616e+01}
& \runid
& 46 
& \num{1.167873e+03}
& \num[scientific-notation=fixed,retain-zero-exponent=false]{6.107816e+01}\%
& \num{1.796794e+01}
\\\bottomrule
\end{tabular}
\end{tiny}
\end{table}
\end{landscape}

\ipoint{Observations} The most important observation is that the nested preconditioner $\mat{P}_{\text{2L}}$ is very effective; it allows us to significantly reduce the number of iterations especially when turning to low regularization parameters. Our new scheme results in a speedup by---on average---one order of magnitude, with a peak performance of more than 20x.

The results in \tabref{t:preconditioners-kktsolve-solutionerror} demonstrate that our schemes all converge to the true solution with an error that is at least at the order of the tolerance used to invert the KKT system, i.e. $\bigO(\num{1E-12})$. The solver does not seem to be sensitive to the CFL number used for the SL method.

The results in \figref{f:preconditioners-kktsolve-convergence-plots} suggest that there are dramatic differences in the performance of our preconditioners; the number of iterations reduces significantly for $\mat{P}_{\text{2L}}$. We can for instance reduce the number of iterations from 310 (\runref{102}) to 7 (\runref{105}). The differences in time to solution, however, are less pronounced. Our spectral discretization makes it in general extremely challenging to design a preconditioner that is more effective than $\mat{P}_{\text{REG}}$ given its ideal application and construction costs; inverting and applying $\mat{P}_{\text{REG}}$ is only at the cost of a spectral diagonal scaling (see \secref{s:precond} for details). The regularization preconditioner is effective for smooth problems and large regularization parameters $\beta_v$ (see the first column and, e.g., \runref{29} (HAND images) or \runref{30} (BRAIN images) in in \tabref{t:preconditioners-kktsolve-convergence}). We can observe that this preconditioner becomes less effective as we decrease $\beta_v$ (see also \cite{Mang:2016a,Mang:2015a}). For example, the number of iterations increases from 33 to 92 to 279 if we reduce $\beta_v$ from \num{1E-1} to \num{1E-2}, and finally to \num{1E-3} (\runref{57}, \runref{66}, and \runref{75} in \tabref{t:preconditioners-kktsolve-convergence}, respectively). We can reduce the number of iterations by more than one order of magnitude if we use the nested preconditioner $\mat{P}_{\text{2L}}$. If we use a PCG method with a tolerance of $\num{1E-7}$ to compute the action of the inverse of $\mat{P}_{\text{2L}}$ the preconditioner is almost ideal, i.e., the number of iterations is independent of $\beta_v$ and the grid size $\vect{n}_x$ (see e.g., \runref{14}, \runref{41}, \runref{68}, and \runref{95} in \tabref{t:preconditioners-kktsolve-convergence}). The low tolerance (in our case one order of magnitude smaller than the tolerance we use to solve the reduced space KKT system) to compute the action of the inverse of $\mat{P}_{\text{2L}}$ results in significant application costs. Despite this increase in application costs we can---already for the present two-dimensional prototype implementation---reduce the time to solution for most of the test problems (see, e.g., \runref{67} or \runref{95} in \tabref{t:preconditioners-kktsolve-convergence}). A significant factor is the SL scheme. We can further reduce the CPU time if we replace the PCG method for computing the action of the inverse of $\mat{P}_{\text{2L}}$ by a CHEB method with a fixed number of iterations. We can see that the effectiveness of this scheme is almost independent of the grid size $\vect{n}_x$ (compare, e.g., \runref{17}, \runref{44}, \runref{71}, and \runref{98} in \tabref{t:preconditioners-kktsolve-convergence}). Also, given that we use a fixed number of iterations the percentage of CPU time spent on applying the preconditioner remains almost constant for $\beta_v$ fixed. The nested preconditioner becomes less effective if we reduce $\beta_v$ from $\num{1E-2}$ to $\num{1E-3}$; the number of iterations increases, which in turn makes the speedup less pronounced (see, e.g., \runref{98} vs. \runref{107} or \runref{99} vs. \runref{108} in \tabref{t:preconditioners-kktsolve-convergence}). Although the speedup varies from case to case, we can see that $\mat{P}_{\text{2L}}$ in combination with CHEB(10) outperforms our original scheme for all experiments.

\ipoint{Conclusions} Our nested preconditioner allows us to reduce the number of iterations by more than one order of magnitude and the time to solution by up to a factor of more than 20. We expect these differences to be more pronounced for an optimized three-dimensional implementation, something we will investigate in a follow up paper.

\subsection{Inverse Solve}

\ipoint{Purpose} To study the rate of convergence of our scheme for the entire inverse solve.

\ipoint{Setup} We consider different test images to study the performance our our numerical scheme (HAND, HEART, BRAIN, UT). We terminate our solver if the gradient is reduced by two orders of magnitude or if the $\ell^\infty$ norm of $\vect{g}^h_k$, $k=1,2,\ldots$, is equal or smaller than $\num{1E-5}$. We consider the regularization preconditioner with an RK2A(0.2) PDE solver and the two level preconditioner with an SL(5) PDE solver. We estimate the eigenvalues for the CHEB method only for the first iteration (zero velocity field). We report results for compressible, near incompressible, and incompressible diffeomorphisms, accounting for different regularization norms ($H^1$-seminorm, $H^2$-seminorm, and $H^3$-seminorm). We study convergence (number of outer iterations and Hessian matvecs) as a function of the grid size, constraints, regularization parameter, and regularization norm. We choose the regularization weights empirically (based on experience from our former work~\cite{Mang:2016a,Mang:2015a}). We report \bipa\item the relative change of the reduced gradient, \item the relative change of the residual between $m_R^h$ and $m_1^h$, \item the number of outer iterations, \item the number of Hessian matvecs, \item the time to solution, and \item the obtained speedup compared to our original scheme.\eipa

\ipoint{Results} We report results for a compressible diffeomorphism ($H^2$-regularization norm) in \tabref{t:preconditioners-invsolve-performance}. We study grid convergence and convergence with respect to different regularization weights $\beta_v$. We report results for an incompressible diffeomorphism in \tabref{t:preconditioners-invsolve-performance-ic} accounting for different regularization norms ($H^1$-seminorm; $H^2$-seminorm; and $H^3$-seminorm). We report results for a near-incompressible diffeomorphism in \tabref{t:preconditioners-invsolve-performance-ric}.

\begin{table}
\caption{Convergence results for the inversion using our formulation for a compressible diffeomorphism ($H^2$-regularization). We report results for registering different sets of images using our original preconditioner ($\mat{P}_{\text{REG}}$; RK2A scheme with a CFL number of 0.2) and the proposed preconditioner ($\mat{P}_{\text{2L}}$; SL scheme with a CFL number of 5 and 10; CHEB method with a fixed number of 10 iterations). We report results for different registration problems: HAND (grid sizes: $128\times128$; $256\times256$; and $512\times512$), HEART (grid size $192\times192$), and BRAIN (grid size: $256\times300$); see \figref{f:registration-problems}. We study convergence as a function of the grid size (HAND images; number of unknowns $n=2n_x^1n_x^2$) and as a function of the regularization parameter $\beta_v$ (HAND, HEART, and BRAIN images). We terminate the inversion if the relative change of the $\ell^\infty$-norm of the reduced gradient is at least two orders of magnitude or if the $\ell^\infty$-norm of the gradient is smaller or equal to $\num{1E-5}$. We report ($i$) the relative change of the reduced gradient $\|\vect{g}^\star\|_{\text{rel}}$, ($ii$) the relative change of the residual $\|\vect{r}\|_{\text{rel}}$ ($L^2$-distance between $m_R$ and $m_1$),  ($iii$) the number of outer iterations,  ($iv$) the number of Hessian matvecs, ($v$) the time to solution, and ($vi$) the speedup compared to our original scheme.
\label{t:preconditioners-invsolve-performance}}
\centering
\resetrunid
\begin{scriptsize}
\begin{tabular}{rlrlrrllrRlL}
\toprule
  $n$ & $\beta_v$ & run & $\mat{P}$ & PDE solver & PC solver & $\|\vect{g}^{\star}\|_{\text{rel}}$ & $\|\vect{r}\|_{\text{rel}}$ & iter & matvecs & time & speedup
\\\midrule
\mcol{12}{HAND}
\\\midrule
  \num{32768}
& \num{1.000000e-01}
& \runid
& $\mat{P}_{\text{REG}}$
& RK2A(0.2)
& ---
& \num{4.849708e-03}
& \num{2.418591e-01}
& 8
& 58
& \num{9.451173e+01}
& ---
\\
&
& \runid
& $\mat{P}_{\text{2L}}$
& SL(5)
& CHEB(10)
& \num{6.817209e-03}
& \num{2.419486e-01}
& 8
& 21
& \num{1.717761e+01}
& \num{5.502030e+00}
\\
& \num{1.000000e-02}
& \runid
& $\mat{P}_{\text{REG}}$
& RK2A(0.2)
& ---
& \num{8.388602e-03}
& \num{1.001290e-01}
& 8
& 97
& \num{1.551698e+02}
& ---
\\
&
& \runid
& $\mat{P}_{\text{2L}}$
& SL(5)
& CHEB(10)
& \num{5.420239e-03}
& \num{9.993405e-02}
& 9
& 30
& \num{2.231333e+01}
& \num{6.954130e+00}
\\
& \num{1.000000e-03}
& \runid
& $\mat{P}_{\text{REG}}$
& RK2A(0.2)
& ---
& \num{8.590485e-03}
& \num{6.478245e-02}
& 11
& 401
& \num{9.887798e+02}
& ---
\\
&
& \runid
& $\mat{P}_{\text{2L}}$
& SL(5)
& CHEB(10)
& \num{8.612365e-03}
& \num{6.500300e-02}
& 11
& 67
& \num{6.630821e+01}
& \num{1.491188e+01}
\\
  \num{131072}
& \num{1.000000e-01}
& \runid
& $\mat{P}_{\text{REG}}$
& RK2A(0.2)
& ---
& \num{9.041345e-03}
& \num{3.324902e-01}
& 12
& 113
& \num{7.712978e+02}
& ---
\\
&
& \runid
& $\mat{P}_{\text{2L}}$
& SL(5)
& CHEB(10)
& \num{7.273864e-03}
& \num{3.320037e-01}
& 13
& 39
& \num{9.313547e+01}
& \num{8.281461e+00}
\\
& \num{1.000000e-02}
& \runid
& $\mat{P}_{\text{REG}}$
& RK2A(0.2)
& ---
& \num{9.668640e-03}
& \num{2.004090e-01}
& 11
& 159
& \num{1.169081e+03}
& ---
\\
&
& \runid
& $\mat{P}_{\text{2L}}$
& SL(5)
& CHEB(10)
& \num{3.520408e-03}
& \num{1.986186e-01}
& 14
& 60
& \num{1.717736e+02}
& \num{6.805941e+00}
\\
& \num{1.000000e-03}
& \runid
& $\mat{P}_{\text{REG}}$
& RK2A(0.2)
& ---
& \num{9.585115e-03}
& \num{1.570570e-01}
& 17
& 758
& \num{1.055508e+04}
& ---
\\
&
& \runid
& $\mat{P}_{\text{2L}}$
& SL(5)
& CHEB(10)
& \num{8.435886e-03}
& \num{1.568643e-01}
& 18
& 150
& \num{6.003528e+02}
& \num{1.758146e+01}
\\
  \num{524288}
& \num{1.000000e-01}
& \runid
& $\mat{P}_{\text{REG}}$
& RK2A(0.2)
& ---
& \num{1.061569e-02}
& \num{3.403908e-01}
& 14
& 134
& \num{6.375926e+03}
& ---
\\
&
& \runid
& $\mat{P}_{\text{2L}}$
& SL(5)
& CHEB(10)
& \num{9.827859e-03}
& \num{3.399367e-01}
& 15
& 46
& \num{4.810214e+02}
& \num{1.325497e+01}
\\
& \num{1.000000e-02}
& \runid
& $\mat{P}_{\text{REG}}$
& RK2A(0.2)
& ---
& \num{1.011627e-02}
& \num{2.108768e-01}
& 13
& 208
& \num{1.219506e+04}
& ---
\\
&
& \runid
& $\mat{P}_{\text{2L}}$
& SL(5)
& CHEB(10)
& \num{1.088959e-02}
& \num{2.113034e-01}
& 16
& 65
& \num{9.723693e+02}
& \num{1.254159e+01}
\\
& \num{1.000000e-03}
& \runid
& $\mat{P}_{\text{REG}}$
& RK2A(0.2)
& ---
& \num{1.113378e-02}
& \num{1.653934e-01}
& 19
& 853
& \num{8.425364e+04}
& ---
\\
&
& \runid
& $\mat{P}_{\text{2L}}$
& SL(5)
& CHEB(10)
& \num{1.127544e-02}
& \num{1.649933e-01}
& 23
& 171
& \num{3.780907e+03}
& \num{2.228398e+01}
\\\midrule

\mcol{12}{HEART}
\\\midrule
  \num{73728}
& \num{1.000000e-02}
& \runid
& $\mat{P}_{\text{REG}}$
& RK2A(0.2)
& ---
& \num{9.100919e-03}
& \num{8.007875e-01}
& 20
& 473
& \num{5.463022e+02}
& ---
\\
&
& \runid
& $\mat{P}_{\text{2L}}$
& SL(5)
& CHEB(10)
& \num{9.487724e-03}
& \num{7.996471e-01}
& 20
& 105
& \num{1.628020e+02}
& \num{3.355623e+00}
\\
& \num{1.000000e-03}
& \runid
& $\mat{P}_{\text{REG}}$
& RK2A(0.2)
& ---
& \num{9.302757e-03}
& \num{5.089111e-01}
& 31
& 1659
& \num{3.771103e+03}
& ---
\\
&
& \runid
& $\mat{P}_{\text{2L}}$
& SL(5)
& CHEB(10)
& \num{9.917084e-03}
& \num{5.089513e-01}
& 31
& 410
& \num{7.518387e+02}
& \num{5.015840e+00}
\\

& \num{1.000000e-04}
& \runid
& $\mat{P}_{\text{REG}}$
& RK2A(0.2)
& ---
& \num{9.807502e-03}
& \num{2.980640e-01}
& 81
& 14455
& \num{6.174050e+04}
& ---
\\
&
& \runid
& $\mat{P}_{\text{2L}}$
& SL(5)
& CHEB(10)
& \num{8.864356e-03}
& \num{2.959864e-01}
& 76
& 2865
& \num{6.449937e+03}
& \num{9.572264e+00}
\\\midrule
\mcol{12}{BRAIN}
\\\midrule
  \num{153600}
& \num{1.000000e-01}
& \runid
& $\mat{P}_{\text{REG}}$
& RK2A(0.2)
& ---
& \num{9.052242e-03}
& \num{4.817233e-01}
& 21
& 269
& \num{1.753250e+03}
& ---
\\
&
& \runid
& $\mat{P}_{\text{2L}}$
& SL(5)
& CHEB(10)
& \num{9.164209e-03}
& \num{4.815864e-01}
& 21
& 70
& \num{1.693690e+02}
& \num{1.035166e+01}
\\
&
& \runid
& $\mat{P}_{\text{2L}}$
& SL(10)
& CHEB(10)
& \num{9.223477e-03}
& \num{4.818267e-01}
& 21
& 70
& \num{1.498461e+02}
& \num{1.170034e+01}
\\
& \num{1.000000e-02}
& \runid
& $\mat{P}_{\text{REG}}$
& RK2A(0.2)
& ---
& \num{8.754648e-03}
& \num{3.210382e-01}
& 74
& 2645
& \num{2.604438e+04}
& ---
\\
&
& \runid
& $\mat{P}_{\text{2L}}$
& SL(5)
& CHEB(10)
& \num{9.292594e-03}
& \num{3.211066e-01}
& 79
& 619
& \num{2.090912e+03}
& \num{1.245599e+01}
\\
&
& \runid
& $\mat{P}_{\text{2L}}$
& SL(10)
& CHEB(10)
& \num{8.989550e-03}
& \num{3.217028e-01}
& 80
& 624
& \num{1.905923e+03}
& \num{1.366497e+01}
\\
& \num{1.000000e-03}
& \runid
& $\mat{P}_{\text{REG}}$
& RK2A(0.2)
& ---
& \num{8.998078e-03}
& \num{2.048439e-01}
& 110
& 10306
& \num{1.208681e+05}
& ---
\\
&
& \runid
& $\mat{P}_{\text{2L}}$
& SL(5)
& CHEB(10)
& \num{9.535745e-03}
& \num{2.050564e-01}
& 74
& 1156
& \num{5.145124e+03}
& \num{2.349178e+01}
\\
&
& \runid
& $\mat{P}_{\text{2L}}$
& SL(10)
& CHEB(10)
& \num{9.214783e-03}
& \num{2.053885e-01}
& 76
& 1239
& \num{4.525267e+03}
& \num{2.670961e+01}
\\\bottomrule
\end{tabular}
\end{scriptsize}
\end{table}

\begin{table}
\caption{Convergence results for the inversion using our formulation for a fully incompressible diffeomorphism (linear Stokes regularization). We report results for registering the UT images (see \figref{f:registration-problems}; grid size: $256\times 256$) using our original preconditioner ($\mat{P}_{\text{REG}}$; RK2A scheme with a CFL number of 0.2) and the proposed preconditioner ($\mat{P}_{\text{2L}}$; SL scheme with a CFL number of 5; CHEB method with a fixed number of 10 iterations). We consider different regularization norms: an $H^1$-seminorm; an $H^2$-seminorm; and an $H^3$-seminorm (from top to bottom). We terminate the inversion if the change in the $\ell^\infty$-norm of the reduced gradient $\vect{g}^h_k$, $k=1,2,\ldots$, is at least two orders of magnitude or if the $\ell^\infty$-norm of $\vect{g}^h_k$ is smaller or equal to $\num{1E-5}$. We report ($i$) the relative change of the reduced gradient $\|\vect{g}^\star\|_{\text{rel}}$, ($ii$) the relative change of the residual $\|\vect{r}\|_{\text{rel}}$ ($L^2$-distance between $m_R$ and $m_1$),  ($iii$) the number of outer iterations,  ($iv$) the number of Hessian matvecs,  ($v$) the time to solution, and  ($vi$) the speedup compared to our original scheme.
\label{t:preconditioners-invsolve-performance-ic}}
\centering
\begin{scriptsize}
\resetrunid
\begin{tabular}{lllrrllrRlL}
\toprule
run & norm & $\mat{P}$ & PDE solver & PC solver & $\|\vect{g}^{\star}\|_{\text{rel}}$ & $\|\vect{r}\|_{\text{rel}}$ & iter & matvecs & time & speedup
\\\midrule
\runid & $H^1$ & $\mat{P}_{\text{REG}}$ & RK2A(0.2) & ---      & \num{8.038571e-03} & \num{1.396907e-02} & 12 & 137 & \num{1.199937e+03} & ---                \\
\runid &       & $\mat{P}_{\text{2L}}$  & SL(5)     & CHEB(10) & \num{6.400500e-03} & \num{1.382715e-02} & 13 &  43 & \num{8.877478e+01} & \num{1.351664e+01} \\
\runid & $H^2$ & $\mat{P}_{\text{REG}}$ & RK2A(0.2) & ---      & \num{9.050688e-03} & \num{1.903253e-01} & 17 & 177 & \num{1.376574e+03} & ---                \\
\runid &       & $\mat{P}_{\text{2L}}$  & SL(5)     & CHEB(10) & \num{8.814947e-03} & \num{1.902670e-01} & 16 &  53 & \num{1.056234e+02} & \num{1.303285e+01} \\
\runid & $H^3$ & $\mat{P}_{\text{REG}}$ & RK2A(0.2) & ---      & \num{9.136639e-03} & \num{6.237208e-01} & 34 & 402 & \num{2.597402e+03} & ---                \\
\runid &       & $\mat{P}_{\text{2L}}$  & SL(5)     & CHEB(10) & \num{9.300021e-03} & \num{6.235015e-01} & 36 & 111 & \num{2.222358e+02} & \num{1.168759e+01}
\\\bottomrule
\end{tabular}
\end{scriptsize}
\end{table}

\begin{table}
\caption{Convergence results for the inversion using our formulation for a near-incompressible diffeomorphism (linear Stokes regularization). We consider the HAND images in \figref{f:registration-problems}. We report results for our original preconditioner ($\mat{P}_{\text{REG}}$; RK2A scheme with a CFL number of 0.2) and the proposed preconditioner ($\mat{P}_{\text{2L}}$; SL scheme with a CFL number of 5; CHEB method with a fixed number of 10 iterations). We terminate the inversion if the change in the $\ell^\infty$-norm of the reduced gradient $\vect{g}^h_k$, $k=1,2,\ldots$, is at least two orders of magnitude or if the $\ell^\infty$-norm of $\vect{g}^h_k$ is smaller or equal to $\num{1E-5}$. We report ($i$) the relative change of the reduced gradient $\|\vect{g}^\star\|_{\text{rel}}$, ($ii$) the relative change of the residual $\|\vect{r}\|_{\text{rel}}$ ($L^2$-distance between $m_R^h$ and $m_1^h$),  ($iii$) the number of outer iterations,  ($iv$) the number of Hessian matvecs,  ($v$) the time to solution, and  ($vi$) the speedup compared to our original scheme.
\label{t:preconditioners-invsolve-performance-ric}}
\centering
\begin{scriptsize}
\resetrunid
\begin{tabular}{llrlrrllrRlL}
\toprule
$\beta_v$ & $\beta_w$ & run & $\mat{P}$ & PDE solver & PC solver & $\|\vect{g}^{\star}\|_{\text{rel}}$ & $\|\vect{r}\|_{\text{rel}}$ & iter & matvecs & time & speedup
\\\midrule
\num{1.000000e-01} & \num{1.000000e-03} & \runid & $\mat{P}_{\text{REG}}$ & RK2A(0.2) & ---      & \num{9.197854e-03} & \num{1.549709e-01} & 10 & 119 & \num{1.759554e+02} & ---                \\
                   &                    & \runid & $\mat{P}_{\text{2L}}$  & SL(5)     & CHEB(10) & \num{5.549841e-03} & \num{1.565259e-01} & 10 &  29 & \num{2.402747e+01} & \num{7.323093e+00} \\
                   & \num{1.000000e-04} & \runid & $\mat{P}_{\text{REG}}$ & RK2A(0.2) & ---      & \num{9.487744e-03} & \num{1.483498e-01} &  9 &  99 & \num{1.633297e+02} & ---                \\
                   &                    & \runid & $\mat{P}_{\text{2L}}$  & SL(5)     & CHEB(10) & \num{9.020463e-03} & \num{1.502411e-01} &  9 &  25 & \num{1.818720e+01} & \num{8.980475e+00} \\
\num{1.000000e-02} & \num{1.000000e-03} & \runid & $\mat{P}_{\text{REG}}$ & RK2A(0.2) & ---      & \num{9.137653e-03} & \num{6.562309e-02} & 14 & 731 & \num{1.504089e+03} & ---                \\
                   &                    & \runid & $\mat{P}_{\text{2L}}$  & SL(5)     & CHEB(10) & \num{9.570499e-03} & \num{6.597113e-02} & 13 &  60 & \num{7.052659e+01} & \num{2.132655e+01} \\
                   & \num{1.000000e-04} & \runid & $\mat{P}_{\text{REG}}$ & RK2A(0.2) & ---      & \num{8.521436e-03} & \num{5.312758e-02} & 13 & 513 & \num{1.091594e+03} & ---                \\
                   &                    & \runid & $\mat{P}_{\text{2L}}$  & SL(5)     & CHEB(10) & \num{9.583670e-03} & \num{5.367633e-02} & 13 &  60 & \num{6.284208e+01} & \num{1.737043e+01}
\\\bottomrule
\end{tabular}
\end{scriptsize}
\end{table}

\ipoint{Observations} The most important observation is that our solver remains effective for the entire inversion irrespective of the regularization weights, norms, and grid size.

The average speedup compared to the stabilized version of our original solver is about 10x (see, e.g., \runref{1} through \runref{6} in \tabref{t:preconditioners-invsolve-performance-ic}) with a peak performance of more than 20x (see, e.g., \runref{17} vs. \runref{18} and \runref{31} vs. \runref{33}  in \tabref{t:preconditioners-invsolve-performance} or \runref{5} vs. \runref{6} in \tabref{t:preconditioners-invsolve-performance-ric}). We can, e.g., reduce the time to solution from $\sim$3 hours to 10 minutes for a $256\times256$ image (\runref{11} vs. \runref{12} in \tabref{t:preconditioners-invsolve-performance}). We can also infer that the reduced accuracy in time does not significantly affect the overall rate of convergence of our solver; the number of outer iterations remains almost constant.

Potential options to further improve our scheme are a re-estimation of the eigenvalues during the solution process once we have made significant progress (i.e., the velocity field changed drastically) and an increased accuracy for the evaluation of the gradient and the objective. That is, we currently use the same accuracy for evaluating the Hessian and the gradient. We might be able to further improve the overall accuracy and maybe convergence if we use a more accurate SL scheme when evaluating the reduced gradient.

\ipoint{Conclusions} Our experiments suggest that our improved solver remains effective irrespective of the regularization norm, regularization weight, or grid size. We can achieve good performance for our compressible, incompressible and near-incompressible formulations for constrained diffeomorphic image registration. We obtain a speedup of about 10x with a peak performance of 20x compared to the stabilized version of our original solver.

\section{Conclusions}
\label{s:conclusions}

With this paper we follow up on our former work on constrained diffeomorphic image registration~\cite{Mang:2016a,Mang:2015a}. We have provided an improved numerical scheme to efficiently solve the registration problem. Our solver features a semi-Lagrangian formulation, which---combined with a two-level preconditioner for the reduced space KKT system---provides a one order of magnitude speedup compared to our original solver~\cite{Mang:2016a,Mang:2015a}.

We have originally described our Newton--Krylov solver in~\cite{Mang:2015a}; this includes a comparison against a first order gradient descent scheme still predominantly used in many diffeomorphic registration algorithms that operate on velocity fields; see, e.g., \cite{Beg:2005a,Hart:2009a,Vialard:2012a}.\footnote{The control equation in \eqref{e:control-elim} corresponds to the reduced $L^2$ gradient, i.e., the variation of the Lagrangian $\F{L}$ in \eqref{e:lagrangian} with respect to $\vect{v}$; we use the gradient in the Sobolev space induced by the regularization operator in our gradient descent scheme in~\cite{Mang:2015a}; see also~\cite{Beg:2005a,Hart:2009a}.} We have extended our original formulation~\cite{Mang:2015a} for constrained diffeomorphic image registration in~\cite{Mang:2016a}. The work in~\cite{Mang:2016a} features a preliminary study of registration quality as a function of regularization norms and weights. The present work focuses on numerical aspects of our solver. Our contributions are:
\begin{itemize}
\item The implementation of an unconditionally stable SL scheme for constrained diffeomorphic image registration.
\item The implementation of a two-level preconditioner for the system in~\eqref{e:newtonsteprs}.
\item A detailed numerical study of our new, improved solver.
\end{itemize}

We perform numerical tests on various synthetic and real-world datasets to study \bipa\item the convergence behavior of our forward solver, \item the adjoint errors of our schemes, \item the grid convergence of our preconditioner, \item the sensitivity of our preconditioner to changes in terms of the regularization parameters, and \item the overall convergence of our solver with respect to different choices for the preconditioner, forward solver, regularization norms, and constraints\eipa. We found that
\begin{itemize}
\item Our original solver (spectral discretization in combination with an RK2 scheme; \cite{Mang:2016a,Mang:2015a}) can become unstable, even for smooth problems (see, e.g., \runref{36} in \tabref{t:selfconvergence-hyperbolic-pdesolvers}).
\item Our new solver (spectral discretization in combination with an RK2A scheme and our SL scheme) remains stable for all considered test cases.
\item The SL scheme results in an order of magnitude speedup due to its unconditional stability compared to the RK2A scheme, subject to a reduction in numerical accuracy (see, e.g., \runref{17} vs. \runref{20} in \tabref{t:selfconvergence-hyperbolic-pdesolvers} or \runref{53} vs. \runref{66} in \tabref{t:convergence-SL-to-RK2A-transporteqs}; we loose two digits accuracy by increasing the CFL number from 0.2 to 5, i.e., by switching from RK2A(0.2) to SL(5)). Our numerical study suggests that this reduction in accuracy is not critical with respect to the overall performance of our Newton--Krylov solver.
\item Our new scheme delivers a reduction in the number of inner iterations (i.e., the solution of the reduced space KKT system) by more than one order of magnitude (e.g., 8 vs. 279 iterations (see \runref{78} vs. \runref{75} in \tabref{t:preconditioners-kktsolve-convergence}) or 7 vs. 310 iterations (see \runref{105} vs. \runref{102} in \tabref{t:preconditioners-kktsolve-convergence}); see also \figref{f:preconditioners-kktsolve-convergence-plots}). More importantly, we observe a speedup of, on average, 10x up to more than 20x (see, e.g., \runref{93} vs. \runref{99} in \tabref{t:preconditioners-kktsolve-convergence} for an individual solve of the KKT system, and \runref{32} vs. \runref{33} in \tabref{t:preconditioners-invsolve-performance}, or \runref{5} vs. \runref{6} in \tabref{t:preconditioners-invsolve-performance-ric} for the entire inversion).
\end{itemize}

Our algorithm can be used in other applications besides medical imaging, such as weather prediction and ocean physics (for tracking Lagrangian tracers in the oceans)~\cite{Kalnay:2002a} or reconstruction of porous media flows~\cite{Fohring:2014a}. Although our method is highly optimized for regular grids with periodic boundary conditions, many aspects of our algorithm carry over. Our current Matlab prototype implementation is not yet competitive with efficient, highly optimized implementations for diffeomorphic image registration in terms of runtime~\cite{Vercauteren:2009a}.\footnote{We provide a more detailed study in~\cite{Mang:2016a}. Here, we show that we can outperform existing approaches for diffeomorphic image registration in terms of registration quality with our new formulation.} Even with the speedup we could achieve here, we are still not competitive with the (highly optimized multi-core) implementation of the algorithm presented in~\cite{Vercauteren:2009a}. We expect this to change for the implementation of our solver for the three-dimensional case (which will feature highly-optimized implementations of the computational kernels of our solver dedicated to multi-core platforms and the design of efficient grid, scale, and parameter continuation schemes to further reduce the time-to-solution). We will extend the study in~\cite{Mang:2016a} by comparing our solver to state-of-the-art implementations of other groups (e.g.,~\cite{Vercauteren:2009a,Vialard:2012a,Avants:2008a}) in terms of time-to-solution, registration quality, and inversion accuracy in the three-dimensional setting, something we are currently actively working on~\cite{Mang:2016c}. We will also investigate other formulations for large deformation diffeomorphic image registration, such as for instance the map based approach in~\cite{Hart:2009a} or the inversion for an initial momentum in~\cite{Vialard:2012a}.

\begin{appendix}

\section{Optimality Conditions}
\label{s:optsys-details}

We can derive the optimality conditions of our problem by computing the first variation of $\F{L}$ with respect to the state, adjoint, and control variables, and applying integration by parts. Our derivation will be formal only. In general, we have to specify the regularity of the underlying objects to ensure existence of an optimal solution. The choices we make for the spaces for the velocity and the images are not independent. We will discuss this in more detail below. If we assume that the objective functional $\F{J}$ and the PDE constraints are continuously differentiable, and satisfy a regularity condition on the constraints (see, e.g.,~\cite{Hinze:2009a}), the following system holds true at a solution $\vect{\phi}^\ast \defeq (m^\ast,\lambda^\ast,p^\ast,w^\ast,\vect{v}^\ast)$ of problem~\eqref{e:ip}:
\begin{subequations}
\label{e:optcond}
\begin{align}
\p_t m^\ast + \igrad m^\ast \cdot \vect{v}^\ast & = 0
&&\text{in}\;\;\Omega \times (0,1],
\label{e:state}
\\
m^\ast &= m_T
&&{\rm in}\;\; \Omega \times\{0\},
\label{e:state-ic}
\\
-\p_t \lambda^\ast - \idiv (\vect{v}^\ast\lambda^\ast) & = 0
&&\text{in}\;\;\Omega \times [0,1),
\label{e:adj}
\\
\lambda^\ast & = m_R  - m^\ast
&&{\rm in}\;\; \Omega \times\{1\},
\label{e:adj-fc}
\\
\idiv\vect{v}^\ast & = w^\ast
&&\text{in}\;\;\Omega,
\label{e:icconstraint}
\\
\beta_v\D{A}[\vect{v}^\ast] + \igrad p^\ast + \vect{b}^\ast & = 0
&&\text{in}\;\;\Omega,
\label{e:control-v}
\\
\beta_w\D{B}[w^\ast] + p^\ast & = 0
&&\text{in}\;\;\Omega,
\label{e:control-w}
\end{align}
\end{subequations}

\noindent with periodic boundary conditions on $\p\Omega$. We refer to~\eqref{e:state} with initial condition~\eqref{e:state-ic}, to~\eqref{e:adj} with final condition~\eqref{e:adj-fc}, and to~\eqref{e:control-v} and~\eqref{e:control-w} as \emph{state} (variation of $\F{L}$ with respect to $\lambda$), \emph{adjoint} (variation of $\F{L}$ with respect to $m$), and \emph{control} (variation of $\F{L}$ with respect to $\vect{v}$ and $w$) equations, respectively. The differential operator $\D{A}$ in~\eqref{e:control-v} corresponds to the first variation of the $H^k$-regularization norms in~\eqref{e:varopt:regv} (see~\eqref{e:diffopregv}). The operator $\D{B}$ in~\eqref{e:control-w} corresponds to the first variation of the regularization operator for $w$.

We can completely eliminate the variables $w$ and $p$, the control equation~\eqref{e:control-w}, and the constraint~\eqref{e:icconstraint} from~\eqref{e:optcond} by simple algebraic manipulations. This is straightforward for $w=0$ (see~\cite{Mang:2015a}) and becomes slightly more involved for a non-zero $w$ (see~\cite{Mang:2016a}). This elimination introduces the pseudo-differential operator $\D{K}$ in~\eqref{e:control-elim}. We refer to our preceding work for more details~\cite{Mang:2015a,Mang:2016a}. Overall, we will arrive at the optimality conditions presented in~\secref{s:optsys}. Computing variations of the weak form of the optimality conditions in~\secref{s:optsys} (which includes the operators $\D{K}$ arising from the elimination of $w$ and $p$) yields the PDE operators for the Newton step.

\begin{remark}
The derivation of the optimality conditions~\eqref{e:optcond} is formal only. We note that the presentation of existence and uniqueness proofs for an optimal solution of~\eqref{e:optcond} are beyond the scope of the present paper. In order for us to ensure the well-posedness of the forward problem, the differentiability of the objective functional and the constraints, and, ultimately, the existence and uniqueness for an optimal solution of the control problem, we have to make sure that the variables in our control formulation meet certain regularity requirements; we have to specify appropriate function spaces for the input images $m_l$, $l\in\{R,T\}$, and the velocity field $\vect{v}$. Several works of other authors have addressed these theoretical requirements in the context of related (optimal control) formulations; see, e.g.,~\cite{Barbu:2016a,Beg:2005a,Borzi:2002a,Chen:2011a,Chen:2011b,Crippa:2007a,DiPerna:1989a,Lee:2010a,Vialard:2009a,Vialard:2012a}. These, e.g., include results for formulations that model images as functions of bounded variation~\cite{Chen:2011a,Vialard:2009a} or functions of Sobolev regularity~\cite{Borzi:2002a,Beg:2005a,Vialard:2012a}, respectively. They consider $H^1$~\cite{Chen:2011b,Crippa:2007a,DiPerna:1989a,Borzi:2002a}, $H^2$~\cite{Beg:2005a,Lee:2010a,Vialard:2012a}, and $H^3$~\cite{Chen:2011a} regularization models for $\vect{v}$, accounting for incompressible~\cite{Chen:2011a,Chen:2011b} or near-incompressible~\cite{Borzi:2002a,Crippa:2007a} velocities. It has also been suggested to stipulate adequate regularity requirements by introducing a diffusion operator into the transport problem~\cite{Barbu:2016a,Kunisch:2011a}.

In our formulation, we model images as compactly supported, smooth functions; we use appropriate mollification and Gaussian smoothing to ensure that we meet these requirements. Numerically, we control the smoothness of the velocity by adjusting the weights for the regularization operator for $\vect{v}$ to ensure that we obtain a diffeomorphic map (up to numerical accuracy). Our experimental results suggest that we stably converge to a local optimal solution using our formulation. However, we note that we are not aware of a theoretical proof that $H^1$-regularity for $\vect{v}$ and its divergence are sufficient to guarantee the existence of an optimal solution of our control problem in the theoretical limit. We also note that we observed instabilities if we stipulate $H^1$-regularity for $\vect{v}$ only without controlling its divergence, in our numerical experiments. Instead of directly controlling the smoothness of $\vect{v}$, we can also (additionally) control the curl of $\vect{v}$; this will add additional regularity to our solution~\cite{Lee:2010a,Amrouche:1998a}. For instance, adding an additional $H^1$-regularization model for the curl of $\vect{v}$ will ensure that $\vect{v}$ is an $H^2$-function. A rigorous proof remains open for future work. Finally, we note that we can change the regularization operators if our formulation for near-incompressible diffeomorphisms does not meet the theoretical requirements; all the derivations and algorithmic features presented here will still apply.
\end{remark}

\section{Stabilized RK2 Scheme}
\label{s:rk2a-derivation}

Here, we present the derivation of the stabilized RK2 scheme introduced in~\secref{s:timeintegration} for the transport equations that appear in our optimality system. We refer to~\cite{Fornberg:1975a,Kreiss:1972a} for a general discussion of this scheme. We start by deriving the antisymmetric form of the forward problem in~\eqref{e:state-pde-elim}. Inserting and subtracting the term $\half{1}\idiv m \vect{v}$ and by using the identity $\idiv m\vect{v} = \igrad m \cdot \vect{v} + m \idiv \vect{v}$ we obtain
\begin{align*}
\p_t m + \half{1}
\big(\igrad m\cdot\vect{v} + \idiv m\vect{v} - m\idiv\vect{v} \big)
& = 0
&&\text{in}\;\;\Omega\times(0,1],
\intertext{with periodic boundary conditions on $\p\Omega$. We use this antisymmetric form as the forward operator. Computing first and second variations using this model results in}
-\p_t \lambda
-\half{1} \big(\idiv\lambda\vect{v}
+ \igrad\lambda\cdot\vect{v}
+ \lambda\idiv\vect{v}\big)
& = 0
&&\text{in}\;\;\Omega\times[0,1),
\\
\p_t \tilde{m} + \half{1}\big(
  \igrad\tilde{m}\cdot\vect{v}
- \tilde{m}\idiv\vect{v}
+ \idiv(\tilde{m}\vect{v}+m\vect{\tilde{v}})
+ \igrad m \cdot\vect{\tilde{v}}
- m\idiv\vect{\tilde{v}}\big)
& = 0
&&\text{in}\;\;\Omega\times(0,1],
\\
-\p_t \tilde{\lambda} - \half{1}\big(
 \igrad\tilde{\lambda}\cdot\vect{v}
+\tilde{\lambda}\idiv\vect{v}
+\idiv(\tilde{\lambda}\vect{v}+\lambda\vect{\tilde{v}})
+\igrad\lambda\cdot\vect{\tilde{v}}
+\lambda\idiv\vect{\tilde{v}}\big)
& = 0
&&\text{in}\;\;\Omega\times[0,1),
\end{align*}

\noindent for the adjoint, incremental state and incremental adjoint equation (notice, that some terms in the above equations will drop for the incompressible case, i.e., for $\idiv\vect{v}=0$). Similarly, we obtain the integro-differential operators
\begin{flalign*}
&&\vect{b}
&= \iut{\half{1}
\big(\lambda\igrad m
- m\igrad\lambda
+ \igrad(\lambda m)\big)}
&&,
\\
\text{and}
&&\vect{\tilde{b}}
&=
\iut{\half{1}\big(
  \lambda\igrad \tilde{m}
- \tilde{m}\igrad\lambda
+ \igrad(\lambda \tilde{m})
+ \tilde{\lambda}\igrad m
- m\igrad\tilde{\lambda}
+ \igrad(\tilde{\lambda} m)
\big)}&&
\end{flalign*}

\noindent for the reduced gradient in \eqref{e:control-elim} and the Hessian matvec in \eqref{e:hessian-matvec}, respectively.

\section{Computational Complexity}
\label{s:computational-complexity}

We report the computational complexity as a function of the number of FFTs we have to compute in \tabref{t:computational-complexity-theoretical}. A comparison of the timings for applying a single FFT and for one cubic spline interpolation step with respect to different grid sizes can be found in \tabref{t:timings-fft-vs-interpolation}. We compute the characteristic $\vect{X}$ for the forward and the adjoint problems only once per iteration. When we evaluate the objective and the gradient we have to solve the state and the adjoint equation. This is when we compute the characteristic; we do not recompute it during the incremental solves. We assign these costs to the solution of the state and adjoint equation.

\begin{table}
\caption
{Computational complexity of our solver for the compressible case. We report this complexity as a function of the number of FFTs and interpolation steps within the key building blocks of our solver. We provide these counts for ($i$) the hyperbolic transport equations that appear in the optimality system (state equation~\eqref{e:state-pde-elim}: SE; adjoint equation~\eqref{e:adj-pde-elim}: AE; incremental state equation~\eqref{e:inc-state-pde-elim}: incSE; incremental adjoint equation~\eqref{e:inc-adj-pde-elim}: incAE), the evaluation of the objective $\F{J}^h$ in~\eqref{e:varopt:objective}, the evaluation of the gradient $\vect{g}^h$ in \eqref{e:control-elim}, and the Hessian matvec in \eqref{e:hessian-matvec}. We report numbers for the full Newton case (FN) and the Gauss--Newton approximation (GN). The costs for evaluating the objective include the costs for the forward solve. We assign the costs of the adjoint solve to the evaluation of the gradient. The costs for the Hessian matvec include the solution of the incSE and incAE.
\label{t:computational-complexity-theoretical}}
\centering
\begin{scriptsize}
\resetrunid
\begin{tabular}{rrrrrrr}
\toprule
 & \mcol{2}{RK2} & \mcol{2}{RK2A} & \mcol{2}{SL} \\\midrule
 & FFTs & IPs & FFTs & IPs & FFTs & IPs
\\\midrule
SE           & $2(d+1)n_t$      & -- & $2(d+1) + 4(d+1)n_t$       & -- & --                    & $d + n_t$          \\
AE           & $2(d+1)n_t$      & -- & $2(d+1) + 4(d+1)n_t$       & -- & $d+1$                 & $d + n_t + 1$      \\
incSE        & $4(d+1)n_t$      & -- & $4(d+1) + 6(d+1)n_t$       & -- & $(d+1)n_t$            & $d+(d+1)n_t$       \\
incAE (FN)   & $2(d+1)n_t$      & -- & $4(d+1) + 6(d+1)n_t$       & -- & $(d+1)n_t$            & $2n_t + 1$         \\
incAE (GN)   & $2(d+1)n_t$      & -- & $2(d+1) + 4(d+1)n_t$       & -- & $d+1$                 & $n_t + 1$          \\
$\D{J}^h$    & $2d + 2(d+1)n_t$ & -- & $2d + 2(d+1) + 4(d+1)n_t$  & -- & $2d$                  & $d + n_t$          \\
$\vect{g}^h$ & $2d + 3(d+1)n_t$ & -- & $2d + 2(d+1) + 7(d+1)n_t$  & -- & $2d + (d+1)(n_t + 1)$ & $d + n_t + 1$      \\
matvec (FN)  & $2d + 8(d+1)n_t$ & -- & $2d + 8(d+1) + 18(d+1)n_t$ & -- & $2d + 4(d+1)n_t$      & $d + (d+3)n_t + 1$ \\
matvec (GN)  & $2d + 7(d+1)n_t$ & -- & $2d + 6(d+1) + 13(d+1)n_t$ & -- & $2d + (d+1)(2n_t+1)$  & $d + (d+2)n_t + 1$
\\\midrule
\end{tabular}
\end{scriptsize}
\end{table}

\begin{table}
\caption{Wall clock times for applying one FFT or one cubic spline interpolation step with respect to different grid sizes. The timings are obtained for the {\tt fftn} and {\tt interp2} functions in Matlab R2013a on a Linux cluster with Intel Xeon X5650 Westmere EP 6-core processors at 2.67GHz with 24GB DDR3-1333 memory.
\label{t:timings-fft-vs-interpolation}}
\centering
\resetrunid
\begin{scriptsize}
\begin{tabular}{rrrR}
\toprule
$n_x^i$ & FFTs               & IPs                & factor             \\\midrule
  16    & \num{3.034000e-05} & \num{4.982930e-03} & \num{1.642363e+02} \\
  32    & \num{3.791000e-05} & \num{2.426170e-03} & \num{6.399815e+01} \\
  64    & \num{1.209200e-04} & \num{2.897990e-03} & \num{2.396618e+01} \\
 128    & \num{3.224100e-04} & \num{5.135130e-03} & \num{1.592733e+01} \\
 256    & \num{8.858000e-04} & \num{2.239586e-02} & \num{2.528320e+01} \\
 512    & \num{3.832400e-03} & \num{1.418876e-01} & \num{3.702317e+01} \\
1024    & \num{1.296024e-02} & \num{6.072565e-01} & \num{4.685534e+01}
\\\bottomrule
\end{tabular}
\end{scriptsize}
\end{table}

\end{appendix}

\end{document}